\documentclass[11pt,reqno]{amsart}
\usepackage[a4paper,margin=30mm]{geometry}

\usepackage{hyperref}

\usepackage{verbatim,upref,amsxtra,amsthm,xcolor,amssymb,amsxtra}

\usepackage[mathscr]{eucal}

\usepackage[initials,nobysame]{amsrefs}

\usepackage{mathtools}
\usepackage[normalem]{ulem}

\usepackage{varioref}

\newcommand{\Z}{\mathbb Z}
\newcommand{\N}{\mathbb N}
\renewcommand{\Re}{\mathsf{Re\, }}

\newcommand{\T}{\mathbb T}
\newcommand{\eps}{\varepsilon}
\newcommand{\bm}{{\boldsymbol m}}
\newcommand{\bn}{{\boldsymbol n}}
\newcommand{\bz}{{\boldsymbol z}}
\newcommand{\bk}{{\boldsymbol k}}
\newcommand{\bl}{{\boldsymbol \ell}}
\newcommand{\bw}{{\boldsymbol w}}
\newcommand{\bo}{{\boldsymbol 0}}

\renewcommand{\S}{{\mathcal S}}

\def\mybf #1{\textbf{\textit{#1}}}

\theoremstyle{plain}
\newtheorem{thm}{Theorem}[section]
\newtheorem{lem}[thm]{Lemma}
\newtheorem{prop}[thm]{Proposition}
\newtheorem{cor}[thm]{Corollary}

\theoremstyle{definition}
\newtheorem{defn}[thm]{Definition}

\newtheorem{exmp}[thm]{Example}
\theoremstyle{remark}
\newtheorem{rem}[thm]{Remark}

\numberwithin{equation}{section}

\date{} 

\begin{document}


\frenchspacing

\setlength\marginparsep{8mm}
\setlength\marginparwidth{20mm}

\title[Bohr chaoticity for principal algebraic actions]{Bohr chaoticity of principal algebraic actions and Riesz product measures}\footnote{This paper has been published online in
{\em  Ergod. Th.} \& {\em Dynam. Sys.}  06 March 2024 doi:10.1017/etds.2024.13}

\author{Ai Hua Fan}
\address{(A. H. Fan) Wuhan Institute for Math \& AI, Wuhan University, Wuhan 430072, China \&
	LAMFA, UMR 7352 CNRS, University of Picardie, 33 rue Saint Leu, 80039 Amiens, France}
\email{ai-hua.fan@u-picardie.fr}

\author{Klaus Schmidt}
\address{(K. Schmidt)
	Mathematics Institute, University of Vienna, Oskar-Morgenstern-Platz 1, 1090 Vienna, Austria}
\email{klaus.schmidt@univie.ac.at}

\author{Evgeny Verbitskiy}
\address{(E. Verbitskiy)
	Mathematical Institute, Leiden University, P.O. Box 9512, 2300 RA Leiden, The Netherlands \& 
Korteweg-de Vries Institute for Mathematics, University of Amsterdam, Postbus 94248, 1090 GE
Amsterdam, The Netherlands  }
\email{evgeny@math.leidenuniv.nl}

	\begin{abstract} 
For a continuous $\mathbb{N}^d$ or $\mathbb{Z}^d$ action on a compact space, we introduce the notion of Bohr chaoticity, which is an invariant of topological conjugacy and which is proved stronger than having positive entropy. We prove that all principal algebraic $\mathbb{Z}$ actions of positive entropy are Bohr-chaotic. The same is proved for principal algebraic actions of $\mathbb{Z}^d$ with positive entropy under the condition of existence of summable homoclinic points.
	\end{abstract}
	
\maketitle

	\section{Introduction}

Bohr chaoticity is a topological invariant introduced in \cite{FFS} for topological dynamical systems. For defining this invariant we recall that a sequence $\boldsymbol w=(w_n)_{n\ge 0}\in \ell^\infty (\mathbb{N},\mathbb{C})$ is a \emph{non-trivial weight sequence} if
	\begin{equation}
	\label{weights}
\limsup_{N\to \infty} \frac{1}{N} \sum_{n=0}^{N-1}|w_n|>0.
	\end{equation}
A non-trivial weight sequence $\boldsymbol{w}$ is \emph{orthogonal} to a topological dynamical system $(X,T)$ if
	\begin{equation}
	\label{orthog}
\lim_{N\to\infty} \frac 1N \sum_{n=0}^{N-1} w_n f(T^nx) =0
	\end{equation}
for every continuous function $f\in C(X)$ and every $x\in X$.

	\begin{defn}[\cite{FFS}]
	\label{d:Bohr chaotic}
A topological dynamical system $(X, T)$ is \emph{Bohr chaotic} if it is non-orthogonal to every non-trivial weight sequence $\boldsymbol{w}\in \ell^\infty(\mathbb{N},\mathbb{C})$. In other words, $(X,T)$ is Bohr chaotic if we can find, for every non-trivial weight sequence $\boldsymbol{w}\in \ell^\infty(\mathbb{N},\mathbb{C})$, a continuous function $g\in C(X)$ and a point $x \in X$ such that
	\begin{equation}
	\label{non-Orth}
\limsup_{N\to \infty} \frac{1}{N} \biggl|\sum_{n=0}^{N-1}w_n g(T^n x)\biggr|>0.
	\end{equation}
	\end{defn}

Bohr chaotic systems must have positive entropy: for example, almost all $(\frac12,\frac12)$ Bernoulli sequences taking values $-1$ and $1$ are orthogonal to every topological dynamical system $(X,T)$ with zero entropy (see \cite{L0}). We list some further basic results on Bohr chaoticity, taken from \cite{FFS}:
	\begin{itemize}
	\item
Any extension of a Bohr chaotic topological dynamical system is Bohr chaotic;
	\item
If a topological dynamical system $(X,T)$ has a nonempty, closed, $T$-invariant subset $Y\subset X$ such that $(Y,T|_Y)$ is Bohr chaotic, then $(X,T)$ is Bohr chaotic;
	\item
No uniquely ergodic dynamical system is Bohr chaotic (this is generalized by Tal \cite{Tal} to systems having at most countably many ergodic measures);
	\item
All affine toral endomorphisms of positive entropy are Bohr chaotic;
	\item
All systems having an $m$-order horseshoe, $m\ge 1$, are Bohr chaotic. By an $m$-order horseshoe $K$ of a system $(X, T)$ we mean a $T^m$-invariant closed non-empty set $K\subset X$ such that the subsystem $(K, T^m)$ is conjugate either to the one-sided shift $(\{0, 1\}^\mathbb{N}, \sigma )$ or to the two-sided shift $(\{0,1\}^\mathbb{Z}, \sigma )$.
	\item
All subshifts of finite type with positive entropy are Bohr chaotic;
	\item
All piecewise monotone $C^1$ interval maps of positive entropy are Bohr chaotic. For example, the $\beta$-shifts with $\beta >1$;
	\item
Every $C^{1+\delta}$ ($\delta>0$) diffeomorphism of a compact smooth manifold admitting an ergodic non-atomic Borel probability invariant measure with non-zero Lyapunov exponents is Bohr chaotic.
	\end{itemize}
The reason for the last two classes is that any such system admits a subsystem which is conjugate to a subshift of finite type of positive entropy (\cite{Katok1980}, \cite{Young1981}).

\smallskip It is interesting to note that for the examples of Bohr chaotic systems constructed in \cite{FFS}, the sets of points $x\in X$ satisfying \eqref{non-Orth} are {large} in the sense that they are of full Hausdorff dimension. Actually, weighted ergodic averages on typical dynamical systems would be multifractal and a study on symbolic spaces is carried out in \cite{Fan2020}.

\medskip In the present paper we extend the notion of Bohr chaoticity from $\Z$- to $\Z^d$-actions and prove that a large class of algebraic dynamical systems --- the so-called principal algebraic actions --- are Bohr chaotic, {provided they have positive entropy}.

\medskip By analogy with \eqref{weights}, we say that a complex sequence $\boldsymbol w=(w_{\boldsymbol n})_{\boldsymbol n \in \mathbb{N}^d}\in \ell^{\infty}(\mathbb{N}^d,\mathbb C)$ is a \emph{non-trivial weight} if
	\begin{displaymath}
\limsup_{N\to \infty} \frac{1}{N^d} \sum_{{\boldsymbol n} \in [0, N-1]^d} |w_{\boldsymbol n}| >0.
	\end{displaymath}
Consider a continuous $\mathbb{N}^d$- or $\mathbb{Z}^d$-action $\alpha \colon \boldsymbol{n}\mapsto \alpha ^{\boldsymbol{n}}$ on some compact space $X$. As in \eqref{orthog} say that a (non-trivial) weight $(w_{\boldsymbol n})_{\boldsymbol n \in \mathbb{N}^d}$ is \emph{orthogonal} to the dynamical system $(X, \alpha )$ if
	\begin{equation}
	\label{eq:chaotic}
\lim_{N\to \infty} \frac{1}{N^d} \sum_{{\boldsymbol n} \in [0, N-1]^d} w_{\boldsymbol n} g (\alpha ^{\boldsymbol n} x) =0
	\end{equation}
for every continuous function $g \in C(X)$ and every points $x\in X$.

	\begin{defn}
	\label{def:Bohrch}
If $\alpha $ is a continuous $\mathbb{N}^d$- or $\Z^d$-action on a compact space $X$ we call $(X, \alpha )$ \emph{Bohr chaotic} if it is not orthogonal to any non-trivial weight, that is to say, if for any non-trivial weight $\boldsymbol w=(w_{\boldsymbol n})_{\boldsymbol n \in \mathbb{N}^d}$ there exist $g \in C(X)$ and $x\in X$ such that
	\begin{equation}
	\label{BC-d}
\limsup_{N\to \infty} \frac{1}{N^d} \biggl|\sum_{{\boldsymbol n} \in [0, N-1]^d} w_{\boldsymbol n} g(\alpha ^{\boldsymbol n} x)\biggr| >0.
	\end{equation}

Note that, if $\alpha $ is a continuous $\mathbb{N}^d$-action on $X$, and if $(\bar{X},\bar{\alpha })$ is the natural extension of $(X,\alpha )$ to a continuous $\mathbb{Z}^d$-action $\bar{\alpha }$ on a compact space $\bar{X}$, then $(\bar{X},\bar{\alpha })$ is Bohr chaotic if and only if the same is true for $(X,\alpha )$. Conversely, if a continuous $\mathbb{Z}^d$-action is Bohr chaotic, it is obviously also Bohr chaotic as an $\mathbb{N}^d$-action.
In view of this last property we focus our attention in much this paper on Bohr chaoticity of $\mathbb{Z}^d$-actions, referring to $\mathbb{N}^d$-actions only where necessary (like in Proposition \ref{thm:0-entropy} or Example \ref{ex:Furstenberg}).
	\end{defn}

As in the $1$-dimensional case one can easily verify the following properties of continuous $\mathbb{Z}^d$-actions $(X,\alpha )$:
	\begin{enumerate}
	\item[(i)]
If $X$ has a closed, $\alpha $-invariant subset $Y$ such that $(Y,\alpha |_Y)$ is Bohr chaotic, then $(X,\alpha )$ is Bohr chaotic;
	\item[(ii)]
	\label{cpe}
If $(X,\alpha )$ has a Bohr chaotic factor $(Y,\beta )$ (i.e., if $(Y,\beta )$ is a Bohr chaotic $\mathbb{Z}^d$-action and there exists a continuous, surjective, equivariant map $\phi \colon X\to Y$), then $(X,\alpha )$ is Bohr chaotic.
	\end{enumerate}
In particular, Bohr chaoticity is an invariant of topological conjugacy.
\medskip

Our main results will be proved by using Riesz product measures borrowed from harmonic analysis and the main technical tool is the notion of $m$-goodness.

\smallskip The paper is organized as follows.
In Section 2, we present algebraic $\Z^d$-actions and their basic properties, state our main results on Bohr chaoticity of principal algebraic $\mathbb{Z}^d$-actions
{{(Theorem \ref{thm:mainZ1} and Theorem \ref{thm:mainZd})}}, and prove that Bohr chaotic algebraic $\mathbb{Z}^d$-actions have to have completely positive entropy (Example \ref{ex:cpe}).
In Section 3, we show that zero-entropy $\mathbb{Z}^d$-actions are not Bohr chaotic. Our main tool, Riesz products, is presented in Section 4 where lacunarity of polynomials is discussed. In Section 5, we prove that any principal algebraic $\Z^d$-action defined by a so-called \emph{$m$-good} polynomial is Bohr chaotic (Theorem \ref{thm:mgood}). Section 6 is devoted to the proof of $m$-goodness for all irreducible polynomials $f\in R_1$ with positive Mahler measure, and Theorem \ref{thm:mainZ1} (for $d=1$) is proved there. Theorem \ref{thm:mainZd} (for $d \ge 2$) is proved in Section 7, where we prove a gap theorem (Theorem \ref{thm:separation}) for irreducible polynomials which admit summable homoclinic points. This gap theorem is of independent interest. In Section 8 we speculate briefly on the necessity of our atorality assumptions for our main results and give some examples of Bohr chaotic principal actions arising from \emph{toral} polynomials.
\medskip

\noindent {\bf Acknowledgement}
The authors would like to thank B. Weiss for valuable discussions and D. Lind for alerting them to an error in the submitted manuscript. The second and the third authors are grateful to Central China Normal University for their hospitality, where part of the work was done. The first author was partly supported by NSF of China (Grant Nos. 11971192 and 12231013).

\section{Algebraic $\Z^d$-actions}
	\label{s:algebraic}
	
	 In this section, we present the principal algebraic $\Z^d$-actions which are our main objects of study and then state our main results
	(Theorem \ref{thm:mainZ1} and Theorem \ref{thm:mainZd}).

An \emph{algebraic $\mathbb{Z}^d$-action} is an action of $\mathbb{Z}^d$ by (continuous) automorphisms of a compact metrizable abelian group. Algebraic $\mathbb{Z}^d$-actions provide a useful source of examples of continuous $\mathbb{Z}^d$-actions with a wide range of properties, both with zero and with positive entropy, and with or without Bohr chaoticity.

We are interested in a particular family of algebraic $\mathbb{Z}^d$-actions, the so-called \emph{cyclic} actions. Denote by $\sigma $ the shift-action of $\mathbb{Z}^d$ on $\mathbb{T}^{\mathbb{Z}^d}$ given by
	\begin{equation}
	\label{eq:shift}
\sigma ^{\boldsymbol{m}}(x)_{\boldsymbol{n}} = x_{\boldsymbol{n}+\boldsymbol{m}}
	\end{equation}
for every $x=(x_{\boldsymbol{n}})_{\boldsymbol{n}\in \mathbb{Z}^d}\in \mathbb{T}^{\mathbb{Z}^d}$. A \emph{cyclic algebraic $\mathbb{Z}^d$-action} is a pair $(X,\alpha _X)$, where $X\subset \mathbb{T}^{\mathbb{Z}^d}$ is a closed, shift-invariant subgroup and $\alpha _X=\sigma |_X$ is the restriction to $X$ of the shift-action $\sigma $ in \eqref{eq:shift}.

In order to describe these actions in more detail we denote by $R_d = \mathbb{Z}[z_1^{\pm1},\dots ,z_d^{\pm1}]$ the ring of Laurent polynomials in the variables $z_1,\dots ,z_d$ with coefficients in $\mathbb{Z}$. Every $f\in R_d$ will be written as $f=\sum_{\bn\in \mathbb{Z}^d}f_\bn {\bz}^\bn$ with $f_\bn\in \mathbb{Z}$ and $\boldsymbol{z}^\bn=z_1^{n_1}\cdots z_d^{n_d}$ for every $\bn=(n_1,\dots ,n_d)\in \mathbb{Z}^d$. The set $\textsf{supp}(f) = \{\bn\in \mathbb{Z}^d\mid f_\bn \ne0\}$ will be called the \emph{support} of $f$, and we set $\|f\|_1=\sum_{\bn\in \mathbb{Z}^d}|f_\bn|$ and $\|f\|_\infty =\max_{\bn\in \mathbb{Z}^d}|f_\bn|$. Following standard terminology we call a nonzero element $f\in R_d$ \emph{primitive} if the greatest common divisor $\gcd(\{f_\bn\mid \bn\in \mathbb{Z}^d\})$ of its coefficients is equal to $1$, and \emph{irreducible} if it is not a product of two non-units in $R_d$. For example, $2\bz^\bn$ is irreducible, but not primitive, while $2(1+\bz^\bn)$ is neither primitive nor irreducible, for every $\bn\in \mathbb{Z}^d$.

\smallskip Every nonzero $f=\sum_{\bn\in \mathbb{Z}^d}f_\bn {\bz}^\bn\in R_d$ defines a surjective group homomorphism $f(\sigma )=\sum_{\boldsymbol{m}\in \mathbb{Z}^d}f_{\boldsymbol{m}}\sigma ^{\boldsymbol{m}}\colon \mathbb{T}^{\mathbb{Z}^d} \to \mathbb{T}^{\mathbb{Z}^d}$. Consider the closed, shift-invariant subgroup
	\begin{equation}
	\label{eq:Xf}
X_f
=\biggl\{ x\in\T^{\mathbb Z^d}\mid \sum_{\bm\in\Z^d} x_{\bn+\bm} f_{\bm} =0 \enspace (\textup{mod}\;1)\enspace \textup{for all}\enspace \bn\in\Z^d\biggr\}
= \ker(f(\sigma ))\subset \mathbb{T}^{\mathbb{Z}^d},
	\end{equation}
and denote by
	\begin{equation}
	\label{eq:alphaf}
\alpha _f = \sigma |_{X_f}
	\end{equation}
the restriction to $X_f$ of the shift-action $\sigma $ on $\mathbb{T}^{\mathbb{Z}^d}$. The dynamical system $(X_f,\alpha _f)$ is called the \emph{principal algebraic action} corresponding to $f\in R_d$.

Formally we can extend this definition of a principal action $(X_f,\alpha _f)$ to include the case $f=0$, the zero polynomial in $R_d$. In this case the definitions \eqref{eq:Xf} -- \eqref{eq:alphaf} reduce to $X_f=\mathbb{T}^{\mathbb{Z}^d}$ and $\alpha _f = \sigma $, i.e., $(X_f,\alpha _f)$ is simply the shift action of $\mathbb{Z}^d$ on $\mathbb{T}^{\mathbb{Z}^d}$.

\smallskip For every cyclic algebraic action $(X,\alpha _X)$ with $X\subsetneq \mathbb{T}^{\mathbb{Z}^d}$, the set
	\begin{equation}
	\label{eq:ideal}
I_X=\{f\in R_d\mid X\subset X_f\}
	\end{equation}
is an ideal in $R_d$ (which is, of course, finitely generated since the ring $R_d$ is Noetherian) and $X= \bigcap_{f \in I_X} X_f$. Conversely, if $I\subset R_d$ is an ideal, generated by $\{f^{(1)},\dots ,f^{(r)}\}$, say, we denote by $(X_I,\alpha _I)$ the cyclic $\mathbb{Z}^d$-action defined by
	\begin{equation}
	\label{eq:XI}
X_I=\bigcap_{f\in I} X_f = \bigcap_{i=1}^r X_{f^{(i)}}\subset \mathbb{T}^{\mathbb{Z}^d}\quad \textup{and}\quad \alpha _I=\sigma |_{X_I},
	\end{equation}
and write $\lambda _I$ for the normalized Haar measure of $X_I$. If the ideal $I\subset R_d$ is principal, $I=(f)$, say, we write $(X_f,\alpha _f)$ instead of $(X_{(f)},\alpha _{(f)})$ and denote by $\lambda _f$ the normalized Haar measure on $X_f$.

	\subsection{Mahler measure}
The topological entropy $h_{\textup{top}}(\alpha _f)$ of a principal algebraic action $(X_f,\alpha _f)$, $f\in R_d\smallsetminus \{0\}$, coincides with its measure-theoretic entropy $h_{\lambda _f}(\alpha _f)$ and is given by the (logarithmic) \emph{Mahler measure} of its defining polynomial $f$:
	\begin{equation}
	\label{eq:mahler}
h_{\textup{top}}(\alpha _f) = \mathsf{m}(f)\coloneqq \int_0^1\cdots \int_0^1 \log |f(e^{2\pi it_1},\dots ,e^{2\pi it_d})|\,dt_1\cdots dt_d.
	\end{equation}

For polynomials in a single variable (i.e., for $f\in R_1$), Mahler measure can be computed using Jensen's formula: let $f(z)=f_0+f_1z +\ldots +f_k z^k$ with $f_0f_k\ne0$ and (complex) roots $\lambda_1,\ldots, \lambda _k$. Then
	\begin{equation}
	\label{eq:mahler1}
\mathsf{m}(f)=\log|f_k| +
\sum_{j:\ |\lambda_j|>1} \log|\lambda_j|.
	\end{equation}
(cf. \cite{LSW}*{p. 597} or \cite{DSAO}*{(16.2)}, or \cite{Young-Jensen}).
The Kronecker lemma \cite{Kronecker}tates that if a polynomial $f\in R_1$ is irreducible, monic and all its roots have absolute value at most $1$,
then $f$ is \emph{cyclotomic}, i.e., for some integer $n$,
	\begin{equation}
	\label{eq:cyclotomic}
f(z)=\Phi_n(z) :=
\prod_{\substack{1\le \ell\le n
	\\
\textsf{gcd}(\ell,n)=1}}
(z-e^{2 \pi i\ell/n}).
	\end{equation}
Since Mahler measure is additive in the sense that
$\mathsf{m}(f\cdot g)=\mathsf{m}(f)+\mathsf{m}(g)$ for all $f,g\in R_d\smallsetminus \{0\}$, any $f\in R_1\smallsetminus \{0\}$ with $\mathsf{m}(f)=0$ must be a product of cyclotomic polynomials.

\medskip A similar statement is true for multivariate polynomials as well:
if $f\in R_d\smallsetminus \{0\}$, then $\mathsf{m}(f)=0$ if and only if $f$ is a product of so-called \emph{generalized cyclotomic polynomials}
	\begin{displaymath}
f(\bz)=\pm {\bz^{\bn_0}}\Phi_{m_1}(\bz^{\bn_1} )
\cdots \Phi_{m_r}(\bz^{\bn_r} )
	\end{displaymath}
for some integers $m_1,\ldots, m_r$ and
$\bn_0,\bn_1,\ldots \bn_r\in \Z^d$ (\cite{Boyd}, \cite{DSAO}*{Theorem 19.5}, \cite{Smyth}).

\medskip We recall the following properties of cyclic algebraic $\mathbb{Z}^d$-action $(X_I,\alpha _I)$ (cf. \cite{DSAO}*{Chapter 6}):
	\begin{itemize}
	\item
The normalized Haar measure $\lambda _{X_I}$ of $X_I$ is shift-invariant;
	\item
If $I\subset R_d$ is nonzero and principal, $I=(f)$, say, the topological entropy of $(X_f,\alpha _f)$ is given by the Mahler measure \eqref{eq:mahler} of $f$; if $I\subset R_d$ contains at least two nonzero elements $f,g$ which are relatively prime to each other (i.e., without a nontrivial common factor), then $h_{\text{top}}(X_I,\alpha _I)=0$;
	\item
If $d>1$, every principal $\mathbb{Z}^d$-action $(X_f,\alpha _f)$ is ergodic (w.r.t. to $\lambda _f$); if $d=1$, a principal $\mathbb{Z}$-action $(X_f,\alpha _f)$ is ergodic if and only if $f$ has no cyclotomic divisor.
	\item
For every nonzero and irreducible $f\in R_d$ the following conditions are equivalent:
	\begin{itemize}
	\item
$\lambda _f$ is mixing under $(X_f,\alpha _f)$,
	\item
$h_{\text{top}}(X_f,\alpha _f)>0$,
	\end{itemize}
	\end{itemize}

\medskip 
\subsection{Main results}
Our main results are the following theorems which will be proved in the Sections \ref{s:d=1} and \ref{s:d>1}.
	\begin{thm}
	\label{thm:mainZ1}
Suppose $f\in R_1\smallsetminus \{0\}$ with $\mathsf{m}(f)>0$. Then
the principal algebraic $\Z$-action $(X_f,\alpha _f)$ is Bohr chaotic.
	\end{thm}

For the higher dimensional case, we need an extra condition:

	\begin {defn}[\cite{LSV}]
	\label{d:atoral}
A nonzero Laurent polynomial $f\in R_d$ is \textit{atoral} if it is not a unit in $R_d$ and its \textit{unitary variety}
	\begin{displaymath}
\mathsf U(f)=\{(t_1,\dots ,t_d)\in \mathbb{T}^d\mid f(e^{2\pi it_1},\dots ,e^{2\pi it_d})=0\}
	\end{displaymath}
of $f$ has dimension $\le d-2$. This includes the possibility that $\mathsf U(f)=\varnothing $, which is equivalent to \emph{expansivity} of the $\mathbb{Z}^d$-action $\alpha _f$. If $\mathsf U(f)$ has dimension $d-1$, $f$ is called \emph{toral}.
	\end{defn}

With this definition, the following is true.

	\begin{thm}
	\label{thm:mainZd}
Suppose that $d\ge2$, and that $f\in R_d$ is atoral. Then $h_\textup{top}(X_f,\alpha _f)>0$ and $(X_f,\alpha _f)$ is Bohr chaotic.
	\end{thm}

We illustrate these definitions with a few examples.

	\begin{exmp}[Toral automorphisms]
	\label{e:toral}
We start with a special case: let $f=f_0+\dots +f_kz^k \in R_1$ with $k\ge 1$ and $f_k=|f_0|=1$. Then the principal $\mathbb{Z}$-action $(X_f,\alpha _f)$ is conjugate to the toral automorphism $(\mathbb{T}^k,A_f)$, where
	\begin{equation}
	\label{eq:companion}
A_f =\left(
	\begin{smallmatrix}
0&1&0&\hdots&0&0&0\vspace{1mm}
	\\
0&0&1&\hdots&0&0&0\vspace{-1mm}
	\\
\vdots&\vdots &\vdots&\ddots&\ddots&\vdots&\vdots\vspace{1mm}
	\\
0&0&0&\hdots&0&1&0\vspace{1mm}
	\\
0&0&0&\hdots&0&0&1\vspace{1mm}
	\\
-f_0&-f_1&-f_2&\hdots&-f_2&-f_{k-2}&-f_{k-1}\vspace{1mm}
	\end{smallmatrix}
\right) \in \textup{GL}(k,\mathbb{Z})
	\end{equation}
is the \emph{companion matrix} of $f$. The map $\phi \colon X_f \to \mathbb{T}^k$, defined by
	\begin{displaymath}
\phi (x)=
	\left(
	\begin{smallmatrix}
x_0\vspace{-1mm}
	\\
\vdots
	\\
x_{k-1}
	\end{smallmatrix}\right)
	\end{displaymath}
for every $x=(x_n)_{n\in \mathbb{Z}}$, implements this conjugacy. We conclude that $(X_f,\alpha _f)$ and thus $(\mathbb{T}^k,A_f)$ is Bohr chaotic if and only if $\mathsf{m}(f)>0$ (cf. \eqref{eq:mahler}).

Consider now an irreducible toral automorphism $T_A\colon \T^k\to\T^k$ defined by a matrix $A\in GL_k(\Z)$. Then the characteristic polynomial $f(z)$ of $A$ is irreducible, and the principal algebraic $\mathbb{Z}$-action $(\T^k,A_f) \cong (X_f,\alpha _f)$ in \eqref{eq:companion} is a finite-to-one factor of $(\T^d,T_A)$. Hence, if $(\T^k, A_f)$ is Bohr chaotic (which is the case if and only if $\mathsf(m)(f)>0$), then $(\T^k,T_A)$ is also Bohr chaotic (as an extension).

If a toral automorphism $T_A\colon \T^k\to\T^k$ with $A\in GL_k(\Z)$ is reducible (i.e., has a proper invariant subtorus $V\subsetneq \mathbb{T}^k$), then the characteristic polynomial $g$ of $T_A|_V$ will be a proper factor of $f$. If $g'$ is one of the irreducible factors of $g$ (and hence of $f$), then the system $(\mathbb{T}^{\textup{deg}(g')},T_{A_{g'}})$ will be Bohr chaotic if and only if $\mathsf{m}(g')>0$, in which case both $(V,T_A|_V)$ and $(\mathbb{T}^d,T_A)$ will be Bohr chaotic. By varying $V$ over the $A$-invariant irreducible subtori of $\mathbb{T}^k$ we see that $(\mathbb{T}^k,T_A)$ is Bohr chaotic if and only if $h_\textup{top}(T_A) = \mathsf{m}(f)>0$.
	\end{exmp}

	\begin{exmp}[Constant polynomials]
	\label{ex:constant}
Suppose that $f=p\in \mathbb{N}$, $p>1$, viewed as a constant polynomial in $R_d$. Then the principal algebraic action $(X_f,\alpha _f)$ arising from this polynomial is the shift-action \eqref{eq:shift} on $X_p\coloneqq \{0,\dots ,p-1/p\}^{\mathbb{Z}^d}$, which is certainly Bohr chaotic. If $p=1$ (or, more generally, if $f$ is a unit in $R_d$), then $X_f$ reduces to a single point and the $\mathbb{Z}^d$-action $\alpha _f$ becomes trivial. By default, $(X_f,\alpha _f)$ is not Bohr chaotic.
	\end{exmp}

	\begin{exmp}[The zero polynomial]
	\label{ex:zero}
So far we have always assumed that the polynomial $f\in R_d$ defining a principal algebraic action $(X_f,\alpha _f)$ is nonzero. If we deviate from this assumption and set $f=0$ (the zero polynomial in $R_d$), then \eqref{eq:Xf} reduces to $X_f=\mathbb{T}^{\mathbb{Z}^d}$, and $\alpha _f$ becomes the shift action $\sigma $ of $\mathbb{Z}^d$ on $\mathbb{T}^{\mathbb{Z}^d}$. For every integer $p>1$, $\mathbb{T}^{\mathbb{Z}^d}$ contains the closed, shift invariant subset $X_p \coloneqq \{0,\dots ,p-1/p\}^{\mathbb{Z}^d}$ in Example \ref{ex:constant}. Since $(X_p,\sigma )$ is Bohr periodic, the same is true for $(X_f,\alpha _f) = (\mathbb{T}^{\mathbb{Z}^d},\sigma )$.
	\end{exmp}

	\begin{rem}
	\label{r:irreducible}
In order to prove our Theorems \ref{thm:mainZ1} and \ref{thm:mainZd} we may assume, without loss of generality, that the polynomial $f$ in either of these theorems is primitive and irreducible.

Indeed, if $f$ is not primitive, then $f = pg$ for some primitive polynomial $g\in R_d$, and $(X_f,\alpha _f)$ has the subsystem $(X_p,\alpha _p= \sigma |_{X_p})$ appearing Example \ref{ex:constant}. Since $(X_p,\alpha _p)$ is Bohr chaotic, the same holds for $(X_f,\alpha _f)$.

Similarly, if $f\in R_d$ is reducible and $\mathsf{m}(f)>0$, then at least one of irreducible factors $g$ of $f$ has positive Mahler measure $\mathsf{m}(g)>0$. If $(X_g,\alpha _g)$ is Bohr chaotic, the Bohr chaoticity of $(X_f,\alpha _f)$ follows immediately from Bohr chaoticity of the subsystem $(X_g,\alpha _g)$. For $d\ge2$ in Theorem \ref{thm:mainZd} we also note that atorality of a polynomial $f\in R_d$ is inherited by all its irreducible factors.
	\end{rem}

	\section{$\Z^d$ actions of zero entropy}
	\label{sec:Bohrchaoticitysec}

Before going to study our Bohr chaotic principal algebraic $\mathbb{Z}^d$-actions, we would like to justify that any continuous $\N^d$- or $\Z^d$-action with zero topological entropy is not Bohr chaotic (Proposition \ref{thm:0-entropy}). This is an immediate consequence of the disjointness completely positive entropy systems and zero entropy systems.
Also we would like to point out that principal $\Z$-actions with zero topological entropy are disjoint from the M\"{o}bius function (Proposition \ref{thm:Mobiusoned}).

\subsection{Zero entropy $\mathbb{Z}^d$-actions are not Bohr chaotic}
	\label{ss:zero}

Consider a measure-pre\-serving $\mathbb{N}^d$- or $\Z^d$-action $\gamma $ on a Lebesgue space $(\Omega , \mu )$, where $\Omega $ is a compact space equipped with its Borel field. We say that the measure-theoretic system $(\Omega , \mu , \gamma )$ has \emph{completely positive entropy} if any non-trivial factor of $(\Omega , \mu , \gamma )$ has positive entropy. Bernoulli systems have complete positive entropy. For $d=1$, the following result is folklore; for $d\ge 1$ we include a proof for completeness, based on a disjointness theorem due to Glasner, Thouvenot and Weiss \cite{GTW2000}*{Theorem 1}.

	\begin{prop}
	\label{thm:0-entropy}
Suppose that $(\Omega, \mu , \gamma)$ has completely positive entropy, $\omega \in \Omega $ is a $\mu $-generic point, and $\phi\in C(\Omega )$ is a continuous function having zero mean. Then $(\phi (\gamma ^{\bn} \omega ))_{\bn \in \mathbb{N}^d}$ is orthogonal to every zero entropy $\mathbb{N}^d$- or $\mathbb{Z}^d$-action $(X, \alpha)$. That is to say, for every $f\in C(X)$ and every $x\in X$, we have
	\begin{equation}
	\label{eq:NBC}
\lim_{N\to\infty}\frac{1}{N^d}\sum_{\bn \in [0, N-1]^d} \phi(\gamma ^{\bn}\omega ) f(\alpha ^\bn x)=0.
	\end{equation}
In particular, continuous $\N^d$- or $\Z^d$-actions with zero topological entropy are not Bohr chaotic.
	\end{prop}

	\begin{proof}
Suppose that for some $f$ and some $x$, there exists a sequence $(N_j)$ tending to infinity such that
	\begin{displaymath}
\ell \coloneqq \lim_{j\to\infty}\frac{1}{N_j^d}\sum_{\bn \in [0, N_j-1]^d} \phi(\gamma ^{\bn}\omega ) f(\alpha ^\bn x) \ne 0.
	\end{displaymath}
We can assume that along this sequence $(N_j)$ the following weak limits of measures exist
	\begin{displaymath}
\lambda \coloneqq \lim_{j\to\infty}\frac{1}{N_j^d}\sum_{\bn \in [0, N_j-1]^d} \delta_{\gamma ^{-\bn}\omega } \negthinspace\times\negthinspace \delta_{\alpha ^{-\bn} x}, \quad \nu \coloneqq \lim_{j\to\infty}\frac{1}{N_j^d}\sum_{\bn \in [0, N_j-1]^d} \delta_{\alpha ^{-\bn} x},
	\end{displaymath}
where $\delta _\omega $ and $\delta _x$ denote the point masses at the points $\omega $ and $x$, respectively. Clearly, the measure $\lambda $ is $\gamma \negthinspace\times\negthinspace \alpha $-invariant, and the projection of $\lambda $ on $X$ is equal to $\nu $. Since $\omega $ is $\mu $-generic, the projection of $\lambda $ onto $\Omega $ is equal to $\mu $. In other words, $\lambda $ is a joining of $\mu $ and $\nu $, where $\nu $ has zero entropy. Since systems of completely positive entropy are disjoint from systems of zero entropy by \cite{GTW2000}*{Theorem 1}, we obtain that $\lambda =\mu \negthinspace\times\negthinspace \nu $. Thus, by the definition of $\lambda $ and the hypothesis that $\mathbb{E}_\mu \phi=0$, we get that
	\begin{displaymath}
\ell = \mathbb{E}_\lambda (\phi \otimes f) = \mathbb{E}_\mu \phi \cdot \mathbb{E}_\nu f =0,
	\end{displaymath}
a contradiction.
	\end{proof}

	\begin{cor}
	\label{ex:cpe}
Let $(X,\alpha )$ be an algebraic $\mathbb{Z}^d$-action which does not have completely positive entropy (w.r.t. the Haar measure $\lambda _X$). Then $(X,\alpha )$ is not Bohr chaotic.
	\end{cor}

	\begin{proof}
If $(X,\alpha )$ does not have completely positive entropy, then \cite{DSAO}*{Theorem 20.8} implies that there exists a nontrivial closed, $\alpha $-invariant subgroup $Y\subset X$ such that the $\mathbb{Z}^d$-action $\alpha _{X/Y}$ induced by $\alpha $ on $X/Y$ has zero entropy. Condition (ii) {{at the beginning of Section \ref{sec:Bohrchaoticitysec}}} ,  combined with Proposition \ref{thm:0-entropy}, shows that $(X,\alpha )$ cannot be Bohr chaotic.
	\end{proof}

	\begin{exmp}[Furstenberg's example]
	\label{ex:Furstenberg}
Let $d=2$, and let $I = (2-z_1, 3-z_2)\subset R_2$. Then $X_I=\{x\in \mathbb{T}^{\mathbb{Z}^2}\mid \sigma ^{(1,0)}x=2x,\,\sigma ^{(0,1)}x=3x\}$, so that $x_{k,l}=2^k3^lx_{(0,0)}$ for every $x\in X_I$ and $(k,l)\in \mathbb{Z}^2$. Since $f^{(2)}=2-z_1$ and $f^{(3)}=3-z_2$ are irreducible and relatively prime to each other, $I$ is a prime ideal, and hence $h_{\textup{top}}(X_I,\alpha _I)=0$ \cite{DSAO}*{Proposition 17.5}.

If $\gamma $ is a continuous $\mathbb{Z}^2$-action on a compact space $\Omega $, $\mu $ is a probability measure on $\Omega $ with completely positive entropy under $\gamma $, $\omega \in \Omega $ is a $\mu $-generic point, and $\phi \in C(\Omega )$ has mean zero, Proposition \ref{thm:0-entropy} shows that
	\begin{displaymath}
\lim_{N\to\infty } \frac{1}{N^2} \sum_{(m,n)\in [0,N-1]^2} \phi (\gamma ^{(m,n)}\omega )h(2^m3^nt)=0
	\end{displaymath}
for every $h\in C(\mathbb{T})$ and $t\in \mathbb{T}$.

\medskip In \cite{Furstenberg1967}, Furstenberg's example was defined as the $\mathbb{N}^2$-action $\alpha $ on $X=\mathbb{T}$ given by
	\begin{displaymath}
\alpha ^{(m,n)}t=2^m3^nt\enspace (\textup{mod\; 1})
	\end{displaymath}
for every $(m,n)\in \mathbb{N}^2$ and $t\in \mathbb{T}$.

We set $\Omega =\mathbb{T}^{\mathbb{N}^d}$, write the coordinates of every $\omega =(\omega _\bn)_{\bn\in \mathbb{N}^d}\in \Omega $ in the form $\omega _\bn = (\omega _\bn^{(1)},\dots ,\omega _\bn^{(d)})$,
and denote by $\gamma $ the one-sided shift-action of $\mathbb{N}^d$ on $\Omega $ (cf. \eqref{eq:shift}). According to Franklin \cite{Franklin1963}, for Lebesgue-\emph{a.e.} $(\beta _1, \dots, \beta _d)$ with $\beta _1>1, \cdots, \beta _d>1$, the point $\beta =(\beta _\bn)_{\bn\in \mathbb{N}^d} \in \Omega $ with $\beta _\bn = (\beta _1^{n_1}\;(\textup{mod}\,1), \dots, \beta _d^{n_d}\;(\textup{mod}\,1))$ for every $\bn\in \mathbb{N}^d$ is Lebesgue-generic for $\gamma $ on $\Omega $. If $\phi \colon \Omega \to \mathbb{C}$ is the map defined by
	\begin{displaymath}
\phi (\omega )= e^{2\pi i(\omega _\mathbf{0}^{(1)}+\dots + \omega _\mathbf{0}^{(d)})},
	\end{displaymath}
then
	\begin{displaymath}
\phi(\gamma ^{\bn}\beta )
=e^{2\pi i(\beta _1^{n_1}+ \dots +\beta_d^{n_d})}
	\end{displaymath}
for every $\bn = (n_1,\dots ,n_d)\in \mathbb{N}^d$. By Proposition \ref{thm:0-entropy}, the sequence $(\phi(\gamma ^{\bn}\beta))_{\bn\in \mathbb{N}^d}$ is almost surely orthogonal to all systems of zero entropy.
Since Furstenberg's example $(\mathbb{T},\alpha )$ described in the preceding paragraph has zero entropy, we obtain the following corollary of Proposition \ref{thm:0-entropy}:

	\begin{cor}
For almost all $(\beta_1, \beta_2)$ with $\beta_1>1$ and $\beta_2>1$,
	\begin{displaymath}
\lim_{N\to\infty}\frac{1}{N^2}\sum_{0\le m, n<N} e^{2\pi i (\beta_1^m + \beta_2^n)} f(2^m 3^n t)=0
	\end{displaymath}
for every continuous function $f \in C(\mathbb{T})$ and every $t\in \mathbb{T}$.
	\end{cor}
	\end{exmp}

	\subsection{M\"obius disjointness and principal actions}
	\label{sec:mobius}

We have just shown that zero entropy $\Z^d$-actions are not Bohr chaotic. In fact, for principal actions the result can be strengthened. We recall that a topological dynamical system $(X,T)$ is \emph{M\"{o}bius disjoint} if
	\begin{equation}
	\label{M-orthog}
\lim_{n\to\infty} \frac 1n \sum_{k=1}^n \mu (k) f(T^kx) =0 \text{\ \ \ for every $f\in C(X)$ and every $x\in X$}.
	\end{equation}

	\begin{prop}
	\label{thm:Mobiusoned}
A zero entropy principal $\Z$-algebraic action $(X_f,\alpha_f)$, $f\in R_1$, is M\"obius disjoint.
	\end{prop}
	\begin{proof}
Since $(X_f,\alpha_f)$ has zero entropy, i.e., $\mathsf{m}(f)=0$,
the Kronecker lemma implies that $f$ has the form
	\begin{equation}
	\label{eq:Kroform}
f(z)= \pm z^{m_0} \Phi_{n_1}(z^{m_1})\cdots \Phi_{n_k}(z^{m_k}),
	\end{equation}
where $m_0\in\Z$, $n_j,m_j\in \N$, $j=1,\ldots,k$, and $\Phi_n$ is the $n$-th cyclotomic polynomial defined in \eqref{eq:cyclotomic}. One immediately concludes from (\ref{eq:Kroform}), that
	$$
f(z)=a_0+a_1z+\ldots+a_{N}z^N\text{ with } |a_0|=|a_N|=1,
	$$
and hence $(X_f,\alpha_f)$ is topologically conjugate to the toral automorphism $(\T^{N},T_A)$, where $T_A:\T^{N}\to \T^{N}$
is a linear automorphism with the matrix $A=A_{f}$ -- the companion matrix of $f$, see Example \ref{e:toral}. However, toral automorphisms with zero entropy are known
to be M\"obius disjoint \cite{Liu}*{Theorem 1.1}. In fact, toral automorphisms, and more generally affine maps of compact abelian groups, are the primary examples motivating Sarnak's (still unproven) conjecture that all topological dynamical systems with zero entropy are M\"{o}bius disjoint (\cite{Sarnak}).
	\end{proof}

\section{Riesz product measures on $X_f$}
	\label{s:riesz}

In this section we start on the proofs of the Theorems \ref{thm:mainZ1} and \ref{thm:mainZd}. As explained in Remark \ref{r:irreducible} we assume from now on --- and without loss in generality --- that the polynomial $f\in R_d$ defining our principal action $(X_f,\alpha _f)$ is primitive, irreducible, and has positive Mahler measure.

For the proofs of these theorems we shall use a class of measures called Riesz products. Firstly, we will recall the general construction of Riesz product measures on arbitrary compact abelian groups. Secondly, we will construct Riesz products on $X_f$ based on lacunary polynomials in the dual group $\widehat{X}_f \subset R_d$.

\subsection{Riesz product measures}
Let $X$ be a compact abelian group with dual group $\widehat X$.

	\begin{defn}[\cite{HZ66}]
	\label{d:dissociate}
An infinite sequence of distinct characters $\Lambda =(\gamma _n)_{n\in\mathbb{N}}=\{\gamma _0,\gamma _1,\linebreak[0]\ldots\,\}\subset \widehat X$ is said to be \mybf{dissociate} if for every $k\ge 1$ and every $k$-tuple $(n_1,n_2,\ldots,n_k)\in \mathbb{N}^k$ of distinct non-negative integers, the equality
	\begin{displaymath}
\gamma _{n_1}^{\eps_1}\gamma _{n_2}^{\eps_2}\ldots \gamma _{n_k}^{\eps_k}=1
	\end{displaymath}
with $\eps_j\in\{-2,-1,0,1,2\}$ for every $j=1,\dots ,k$, implies that
	\begin{displaymath}
\gamma _{n_1}^{\eps_1}=\gamma _{n_2}^{\eps_2}=\ldots =\gamma _{n_k}^{\eps_k}=1.
	\end{displaymath}
Equivalently, $\Lambda $ is dissociate if any character in $\widehat X$ can be represented in at most one way as a finite product $ \gamma _{n_1}^{\eps_1}\gamma _{n_2}^{\eps_2}\ldots \gamma _{n_k}^{\eps_k} $ of elements of $\Lambda $, where all $n_j$ are distinct and $\eps_j\in\{-1, 0,1\}$.
	\end{defn}

Using dissociate sequences of characters, Hewitt and Zuckermann \cite{HZ66} proposed a construction of interesting probability measures -- the so-called Riesz products, generalizing Riesz products on $\T$ constructed by F. Riesz \cite{Riesz} in 1918. More precisely, denote by $\lambda _X$ the Haar measure on $X$. Suppose that
	\begin{itemize}
	\item[(i)]
$\Lambda = (\gamma _n)_{n\ge 0}$ is a dissociate sequence of characters in $\widehat X$,
	\item[(ii)]
$a = (a_n)_{n\ge 0}$ is a sequence of complex numbers such that $ |a_n|\le 1$ for all $n$.
	\end{itemize}
For any $N\ge 0$, denote by $\mu _a^{(N)}$ the measure on $X$ which is absolutely continuous with respect to $\lambda _X$ with density
	\begin{displaymath}
\frac {d\mu _a^{(N)}}{d\lambda _X}(x) = \prod_{n=0}^N\bigl(1 +\Re a_n\gamma _n(x) \bigr).
	\end{displaymath}
It is not very difficult to show that the sequence of measures $(\mu _{a}^{(N)})_{N\ge 0}$ converges weakly; the limiting measure $\mu _a=\lim_N \mu _a^{(N)}$ is called the \emph{Riesz product}, and we denote it as
	\begin{equation}
	\label{eq:notatRiesz}
\mu _a=\prod_{n=0}^\infty \bigl(1 +\Re a_n\gamma _n(x) \bigr).
	\end{equation}
The Riesz product $\mu _a$ is absolutely continuous with respect to the Haar measure $\lambda _X$ if and only if $\sum_{n} |a_n|^2<\infty$
and it is singular to the Haar measure $\lambda _X$ if and only if $\sum_{n} |a_n|^2=\infty$ (see \cite{Peyriere1975}, \cite{Zygmund2002}). We will omit dependence of $\mu _a$ on the sequence $\Lambda $, since $\Lambda $ will usually be fixed.

Since
	\begin{displaymath}
1 +\Re a_n\gamma _n(x) =1 + \frac {a_n}2\gamma _n(x)+ \frac {\overline a_n}2\gamma _n^{-1}(x),
	\end{displaymath}
the Riesz product $\mu _a$, associated to the sequences $\Lambda $ and $a$, can be characterized by the Fourier coefficients $ \widehat{\mu }_a(\gamma )=\int \overline{\gamma }(x) d\mu _a(x)$, $\gamma \in\widehat{X}$, as follows:
	\begin{itemize}
	\item[(a)]
For any finite set of distinct characters $\{\gamma _{n_1},\gamma _{n_2},\ldots, \gamma _{n_k}\}\subset\Lambda $ and any $(\eps_1, \eps_2, \ldots,\linebreak[0]\eps_k)\in \{-1,0,1\}^k$,
	\begin{equation}
	\label{eqn:Riesz}
{\widehat\mu }_a(\gamma _{n_1}^{\eps_1} \gamma _{n_2}^{\eps_2}\cdots \gamma _{n_k}^{\eps_k}) = a_{n_1}^{(\eps_1)} a_{n_2}^{(\eps_2)} \cdots a_{n_k}^{(\eps_k)},
	\end{equation}
where $a_n^{(\eps)}=\frac{a_n}{2}, 0,$ or $\frac{\overline a_n}{2}$, {
depending on whether} $\varepsilon =1,0$, or $-1$;

	\item[(b)]
For any character $\gamma \in\widehat{X}$ not of the form $\gamma _{n_1}^{\eps_1}\gamma _{n_2}^{\eps_2}\cdots \gamma _{n_k}^{\eps_k}$ with $\eps_1,\eps_2,\linebreak[0]\ldots,\eps_k\in \{-1,0,1\}$ as in case (a) above, one has
	\begin{equation}
	\label{eqn:Riesznull}
\widehat{\mu }_a(\gamma )=0.
	\end{equation}
	\end{itemize}

For any two Riesz products $\mu _a$ and $\mu _b$, it is proved in \cite{Peyriere1975} that $\mu _a$ and $\mu _b$ are mutually singular if $\sum |a_n-b_n|^2=\infty$, and mutually equivalent if $\sum |a_n-b_n|^2<\infty$ and $\sup_n|a_n|<1$. For any Riesz product $\mu _a$, it is proved in \cite{Fan1993} that the orthogonal series $\sum c_n (\gamma _n(x)- a_n/2)$ (with $c_n\in \mathbb{C}$) converges $\mu _a$-a.e. if and only if $\sum|c_n|^2<\infty$. Such convergence results will be useful to us in the proofs of Theorem \ref{thm:mainZ1} and Theorem \ref{thm:mainZd}. Riesz products on $\mathbb{T}$ and some generalized Riesz products appear as spectral measures of some dynamical systems (see \cites{Bourgain1993,Ledrappier1970,Queffelec1987}). Riesz products
are tools in harmonic analysis (see \cites{Kahane1970,Katznelson2004,Zygmund2002}).

\subsection{The dual group $\widehat{X}_f$} Before constructing Riesz products on $X_f$,
let us describe the dual group of $X_f$ (cf. \cites{DSAO,LS}). Every Laurent polynomial with integer coefficients
	\begin{displaymath}
h(\bz) =\sum_{\bm \in \mathbb{Z}^d} h_\bm \bz^\bm \in R_d,
	\end{displaymath}
defines a character $\gamma ^{(h)} \in \widehat{\mathbb{T}^{\mathbb{Z}^d}}$, given by
	\begin{displaymath}
\gamma ^{(h)}(x) \coloneqq e^{2\pi i\langle h,x\rangle},
	\end{displaymath}
where
	\begin{displaymath}
\langle h,x\rangle=\sum_{\bm\in\Z} h_\bm x_\bm
	\end{displaymath}
for every $x\in \mathbb{T}^{\mathbb{Z}^d}$. Conversely, every character of $\mathbb{T}^{\mathbb{Z}^d}$ is of the form $\gamma =\gamma ^{(h)}$ for some $h\in R_d$, so that we may identify $\widehat{\mathbb{T}^{\mathbb{Z}^d}}$ with $R_d$. Note, however, that the group operation in $R_d$ is addition, whereas in $\widehat{\mathbb{T}^{\mathbb{Z}^d}}$ it is multiplication:
	\begin{displaymath}
\gamma ^{(h+h')} = \gamma ^{(h)}\gamma ^{(h')}
	\end{displaymath}
for all $h,h'\in R_d$.

\smallskip Since $X_f$ is a subgroup of $\mathbb{T}^{\mathbb{Z}^d}$, every character $\gamma ^{(h)}\in \widehat{\mathbb{T}^{\mathbb{Z}^d}}$, $h\in R_d$, restricts to a character $\tilde{\gamma }^{(h)}\in \widehat{X}_f$. From the definition of $X_f$ in \eqref{eq:Xf} it is clear that, for any two polynomials $h,h'\in R_d$, $\tilde{\gamma }^{(h)} = \tilde{\gamma }^{(h')}$ if and only if $h-h'$ is a multiple of $f$. This allows us to identify the dual group $\widehat{X}_f$ with $R_d/(f)$, where $(f)=R_d\cdot f$ is the principal ideal in $R_d$ generated by $f$:
	\begin{displaymath}
\widehat{X}_f = R_d/(f).
	\end{displaymath}
More generally, if $I\subset R_d$ is an ideal and $X_I$ is given by \eqref{eq:XI}, then
	\begin{displaymath}
\widehat{X}_I= R_d/I.
	\end{displaymath}

\subsection{Lacunary polynomials}

For the construction of Riesz product measures on $X_f$ we have to take a closer look at dissociate families $\Lambda \subset \widehat{X}_f$ in the sense of Definition \ref{d:dissociate}.

	\begin{defn}
	\label{d:m-good}
Given an integer $m\in\mathbb N$, we say that a primitive irreducible polynomial $f\in R_d$ is \emph{$m$-good} if the following conditions hold:
	\begin{itemize}
	\item[(C1)]
The collection of characters
	\begin{displaymath}
\bigl\{ \tilde{\gamma }^{(\boldsymbol{z}^{m\bn})} \mid \bn\in\mathbb{N}^d \bigr\}\subset \widehat{X}_f
	\end{displaymath}
is dissociate. Explicitly, this means that any nonzero polynomial
of the form $g(\bz^m)$ where
	\begin{displaymath}
g(\bz)=\sum_{\bn \in \mathbb{Z}^d} \varepsilon _\bn \bz^{\bn}
	\end{displaymath}
with $\varepsilon _\bn \in \{-2,-1,0,1,2\}$
is not divisible by $f$.
	\item[(C2)]
For any $\bk \in \mathbb{N}^d/ m\mathbb{N}^d$, any two points $\bn \ne \bn'$ in $\mathbb{Z}^d$, and any nonzero polynomial of the form $g(\bz)\coloneqq\sum_{\bn \in \mathbb{Z}^d} \varepsilon _\bn \bz^{\bn}$ with $\varepsilon _\bn \in \{-1,0,1\}$,
the polynomial
	\begin{displaymath}
{\boldsymbol z}^{m\bn+\bk}-{\boldsymbol z}^{m\bn'+\bk} + g(\bz^{m})
	\end{displaymath}
is not divisible by $f$.
	\end{itemize}
	\end{defn}

For a given principal algebraic action $(X_f,\alpha _f)$, where $f$ is $m$-good, Riesz product measures $\mu _a$ can be constructed using the countable dissociate collection of characters $\Lambda =\bigl\{ \tilde{\gamma }^{({\boldsymbol z}^{m\bn})}\mid \ \bn\in\mathbb{N}^d\bigr\}$ (cf. (C1)). The second condition (C2) ensures that any shifted family of characters $\Lambda _{\bk}=\bigl\{ \tilde{\gamma }^{({\boldsymbol z}^{m\bn+\bk})}\mid \ \bn\in\mathbb{N}^d \bigr\}$ (with $\bk\in[0,m-1]^d\smallsetminus\{\mathbf{0}\}$ being fixed) is a $\mu _a$-orthogonal system, as a direct consequence of \eqref{eqn:Riesz} applied with $k=2$ --- a useful property which will help us control the behavior of weighted ergodic averages. As we will see, the coefficient sequence $a$ will be chosen depending on the non-trivial weight sequence $\boldsymbol w$.

\section{$(X_f, \alpha _f)$ is Bohr chaotic when $f$ is $m$-good}

The following theorem will allow us to reduce the proof of Bohr chaoticity of $(X_f, \alpha _f)$ to checking the $m$-goodness of the polynomial $f$.

	\begin{thm}
	\label{thm:mgood}
If a primitive irreducible polynomial $f\in R_d$ with positive Mahler measure is $m$-good, i.e., if the conditions \textup{(C1)} and \textup{(C2)} hold for some positive integer $m$, then $(X_f,\alpha _f)$ is Bohr chaotic.
	\end{thm}

We begin with a simple auxiliary lemma.

	\begin{lem}
	\label{lem:aux}
Let $\alpha $ be a continuous $\Z^d$-action on a compact metrizable space $X$, and let $\bw = (w_\bn)_{\bn \in \mathbb{N}^d}$ be a non-trivial weight. Then $(X,\alpha )$ is not disjoint from $\bw=(w_\bn)$ if and only if for any $\bk\in\mathbb{N}^d$, $(X,\alpha )$ is not disjoint from the weight $\widetilde{\bw}=(\widetilde{w}_\bn)$ defined by $\widetilde{w}_\bn=w_{\bn+\bk}$ for all $\bn\in\mathbb{N}^d$.
	\end{lem}
	\begin{proof}
Introduce the following notation: for a continuous function $\phi$ on $X$ let
	\begin{displaymath}
S_N^{\bw}\phi(x)=\sum_{{\boldsymbol n} \in [0, N-1]^d} w_{\boldsymbol n} \phi(\alpha ^{\boldsymbol n} x).
	\end{displaymath}
For any $\bk\in\mathbb{N}^d$ and for any $x\in X$, 
one has

$$\aligned
|S_{N+\|\bk\|_\infty}^{\bw}\phi(x)-S^{\widetilde{\bw}}_N\phi(\alpha ^{\bk}x)|
&=\left|\sum_{\bn \in [0, N+\|\bk\|_\infty-1]^d} w_{\bn}\phi(\alpha ^{\bn}x)-\sum_{\bn \in [0, N-1]^d} {w}_{\bn+\bk}\phi(\alpha ^{\bn+\bk}x)\right|\\
&\le \|\phi\|_{\infty}\cdot \bigl| [0,N+\|\bk\|_\infty-1]^d\triangle (\bk+[0,N-1]^d)\bigr|.
\endaligned
$$
To finish the proof, it suffices to notice that the cardinality of
the symmetric difference is of order $O(N^{d-1})$.
\end{proof}

	\begin{proof}[Proof of Theorem \ref{thm:mgood}]
Fix $m\in\N$ such that the conditions (C1) and (C2) hold. Assume that $\bw$ is a non-trivial weight (cf. \eqref{BC-d}). Then for some $\bk\in[0,\ldots,m-1]^{d}$, one has
	\begin{equation}
	\label{eq:subseq}
\limsup_{N\to \infty} \frac{1}{N^d} \sum_{{\bn\mid m\bn+\bk} \in [0, N-1]^d} \left|w_{m\bn+\bk}\right| >0.
	\end{equation}
Without loss of generality we can assume $\bk=\boldsymbol 0$. Otherwise, consider the shifted weight $\widetilde{\bw}=(\widetilde{w}_\bn)$ with $\widetilde{w}_\bn=w_{\bn+\bk}$. By Lemma \ref{lem:aux}, $(X_f,\alpha _f)$ is not disjoint from $\bw$ if and only if $(X_f,\alpha _f)$ is not disjoint from $\widetilde\bw$. Thus it sufficient to consider the weight $\widetilde{\bw}$ for which we can assume that \eqref{eq:subseq} holds with $\bk=0$. In the following we consider an arbitrary such weight.

\bigskip \mybf{Step 1. Choice of the function $\phi$ and the point $x$.}\ We are going to show that \eqref{BC-d} holds for $\phi(x)\coloneqq e^{2\pi i x_0}= e^{2\pi i \langle 1,x\rangle}$ and for almost all $x\in X_f$ with respect to an appropriately chosen Riesz product measure. Note that for all $\bn\in\mathbb{N}^d$,
	\begin{displaymath}
\phi(\alpha_f ^{\bn}x) = e^{2\pi i x_{\bn}}= e^{2\pi i \langle \bz^{\bn}, x\rangle}=\gamma ^{(\bz^{\bn})} (x).
	\end{displaymath}

\bigskip \mybf{Step 2. Choice of the measure.}\ By condition (C1), the collection of characters
	\begin{displaymath}
\Lambda \coloneqq\bigl\{ \gamma _{\bn}=\gamma ^{(\bz^{m\bn})}\mid \bn\in\mathbb{N}^d\bigr\}.
	\end{displaymath}
is dissociate. Consider now the following collection of coefficients
	\begin{displaymath}
a \coloneqq\bigl\{ a_{\bn}=e^{-i\,\textup{arg}\, w_{m\bn}}\mid \bn\in\mathbb{N}^d\bigr\}.
	\end{displaymath}
Since $|a_{\bn}|=1$ for all $\bn$, the Riesz product $\mu _a$ in \eqref{eq:notatRiesz} is well defined.

\bigskip \mybf{Step 3. Orthonormality.}
\ For each $\bk\in[0,m-1]^d\smallsetminus \{\boldsymbol 0\}$, consider the following collection of functions
	\begin{displaymath}
\mathcal F_{\bk} \coloneqq\bigl\{ \gamma ^{(\bz^{m\bn+\bk})}(x)=\phi\circ\alpha _f^{m\bn+k}(x)\mid \bn\in\mathbb{N}^d\bigr\}.
	\end{displaymath}
We claim that for each $\bk\in[0,m-1]^d\smallsetminus \{\boldsymbol 0\}$, $\mathcal F_{\bk}$ is orthonormal in $L^2(X_f,\mu _a)$. Indeed, for each $\bn\ne\bn'$, the condition (C2) means that the character corresponding to the polynomial $\bz^{m\bn+\bk}-\bz^{m\bn'+\bk}$:
	\begin{displaymath}
\gamma ^{(\bz^{m\bn+\bk}-\bz^{m\bn'+\bk})}(x)=\gamma ^{(\bz^{m\bn+\bk})}(x) \overline{\gamma ^{(\bz^{m\bn'+\bk})}(x)}
	\end{displaymath}
cannot be expressed as a product of characters in $\Lambda $, and hence using expression \eqref{eqn:Riesznull} for the Fourier coefficients of Riesz products one gets that
	\begin{displaymath}
\int_{X_f} \gamma ^{(\bz^{m\bn+\bk})}(x) \overline{\gamma ^{(\bz^{m\bn'+\bk})}(x)} d\mu _a(x) = \widehat{\mu _a}\bigl(\overline{\gamma ^{(\bz^{m\bn+\bk}-\bz^{m\bn'+\bk})}}\bigr)=0.
	\end{displaymath}
Since $|\gamma ^{(\bz^{m\bn+\bk})}(x)|^2=1$ for all $x$, the orthonormality of $\mathcal F_{\bk}$ is thus proved.

For $\bk=0$, we set
	\begin{displaymath}
\mathcal F_{\boldsymbol 0} \coloneqq\bigl\{ \gamma ^{(\bz^{m\bn})}(x)-\frac{a_{\bn}}2\mid \bn\in\mathbb{N}^d\bigr\}.
	\end{displaymath}
Direct application of formulae \eqref{eqn:Riesz} and \eqref{eqn:Riesznull} immediately gives that the collection of functions $\mathcal F_{\boldsymbol 0}$ is orthogonal in $L^2(X_f,\mu _a)$, and that
	\begin{displaymath}
\int_{X_f} |\gamma ^{(\bz^{m\bn})}(x)|^2d\mu _a(x)=1-\frac{|a_{\bn}|^2}{4}=\frac {3}{4} \,\text{ for all }\,\bn\in\mathbb{N}^d.
	\end{displaymath}

\bigskip \mybf{Step 4. Almost everywhere convergence.}\ Write
	\begin{displaymath}
S_N^{\bw}\phi(x)= \sum_{\bn\in [0,N-1]^d} w_{\boldsymbol n} \phi(\alpha _f^{\boldsymbol n} x) = \sum_{\bk\in [0,m-1]^d}S_{N,\bk}^{\bw}\phi(x),
	\end{displaymath}
where
	\begin{displaymath}
S_{N,\bk}^{\bw}\phi(x)\coloneqq \sum_{\{\bn\mid m\bn+\bk\in [0,N-1]^d\}} w_{m\boldsymbol n+\bk} \phi(\alpha _f^{m\boldsymbol n+\bk}x).
	\end{displaymath}
We claim that for any $\bk\in [0,m-1]^d\smallsetminus\{\boldsymbol 0\}$, one has
	\begin{equation}
	\label{eq:convk}
\frac 1{N^d} S_{N,\bk}^{\bw}\phi(x)\to 0\quad\mu _a\!-\!a.e.,
	\end{equation}
and for $\bk=0$, one has
	\begin{equation}
	\label{eq:convnull}
\frac 1{N^d}\biggl( S_{N,{\boldsymbol 0}}^{\bw}\phi(x)-\frac 12\sum_{\{\bn\mid m\bn\in[0,N-1]^d\}}|w_\bn|\biggr)\to0 \quad\mu _a\!-\!a.e.
	\end{equation}
Now we write
	\begin{displaymath}
\frac {1}{N^d}S_N^{\bw}\phi(x)=\frac {1}{N^d}\biggl(S_N^{\bw}\phi(x) -\frac 12\sum_{\{\bn\mid m\bn\in[0,N-1]^d\}}|w_\bn|\biggr)+\frac {1}{2N^d}\sum_{\{\bn\mid m\bn\in[0,N-1]^d\}}|w_\bn|.
	\end{displaymath}
If \eqref{eq:convk} and \eqref{eq:convnull} are indeed true, the first term in the brackets on the right hand side converges to $0$ for $\mu _a$-almost all $x\in X_f$, and the second term does not converge to $0$ by \eqref{eq:subseq}. Hence, we will be able to conclude that
	\begin{displaymath}
\limsup_{N\to\infty} \frac 1{N^d} \left| S_N^{\bw}\phi(x)\right| >0, \quad \mu _a\!-\!a.e.,
	\end{displaymath}
and thus, that $(X_f,\alpha _f)$ is Bohr chaotic.

Finally, to establish \eqref{eq:convk} and \eqref{eq:convnull}, we will use the following multivariate generalization of the result of Davenport, Erd\"os, and LeVeque \cite{DEL} due to Fan, Fan, and Qiu \cite{FFQ}*{Theorem 6.1}: Suppose that $\{\xi_\bl\mid \bl\in\N^d\}$ is a collection of measurable complex valued uniformly bounded functions on a probability space $(\Omega ,\mathbb{P})$ such that
	\begin{equation}
	\label{eq:condFFQ}
\sum_{N=1}^\infty \frac 1{N}\int_{\Omega } |Z_{N}|^2 d\mathbb{P}<\infty,
	\end{equation}
where
	\begin{displaymath}
Z_N=\frac {1}{N^d}\sum_{\bl\in[0,N-1]^d} \xi_{\bl} \quad (N\ge 1).
	\end{displaymath}
Then $Z_N\to 0$ as $N\to\infty$ $\mathbb{P}$-a.e. on $\Omega $.

\smallskip In particular, if $\{\xi_\bl\mid \bl\in\N^d\}$ are uniformly bounded and orthogonal in $L^2(\Omega ,\mathbb{P})$, then
	\begin{displaymath}
\frac 1{N}\int_{\Omega } |Z_{N}|^2 d\mathbb{P} =\frac 1{N^{2d+1}}\sum_{\bl\in[0,N-1]^d} \int_{X} |\xi_{\bl}|^2 d\mathbb{P} \le \frac {C}{N^{d+1}},
	\end{displaymath}
and hence \eqref{eq:condFFQ} holds for any $d\ge 1$.

\smallskip If we now apply this result to the orthogonal families of bounded functions
	\begin{displaymath}
\mathcal F_{\bk}^{\boldsymbol w}=\Bigl\{ w_{m\bn+\bk}\psi\circ\alpha ^{m\bn+\bk}(x)\mid \bn\in\Z^d_+ \Bigr\},\quad \bk\in[0,m-1]^d\smallsetminus\{\boldsymbol 0\},
	\end{displaymath}
and
	\begin{displaymath}
\mathcal F_{\boldsymbol 0}^{\boldsymbol w}=\Bigl\{ w_{m\bn}\Bigl(\psi\circ\alpha ^{m\bn}(x)-\frac {a_{\bn}}{2}\Bigr) \mid \bn\in\N^d \Bigr\},
	\end{displaymath}
we obtain \eqref{eq:convk} and \eqref{eq:convnull}, and hence, we complete the proof.
	\end{proof}

\section{Bohr chaoticity of $(X_f, \alpha _f)$: the case of $d=1$}
	\label{s:d=1}

In this section we complete the proof of Theorem \ref{thm:mainZ1} in the case where $d=1$: every principal algebraic $\Z$-action $(X_f,\alpha _f)$ with positive entropy is Bohr chaotic. Theorem \ref{thm:mainZ1} will follow from Remark \ref{r:irreducible}, Theorem \ref{thm:mgood} and the following result.

	\begin{thm}
	\label{thm:mgood_d=1}
Every primitive irreducible polynomial $f\in R_1$ with $\mathsf{m}(f)>0$ is $m$-good for some positive integer $m$.
	\end{thm}

The proof of Theorem \ref{thm:mgood_d=1} consists of the following three lemmas.

	\begin{lem}[Preliminary lemma]
	\label{lem:lind}
Let $f(z) = f_0 + f_1 z + \cdots + f_r z^r \in \mathbb{Z}[z]$ be an irreducible polynomial with $r\ge1$ and $f_0f_r\ne 0$. If $\mathsf{m}(f)>0$ then at least one of the following statements is true.
	\begin{enumerate}
	\item
There exists a root of $f$ in $\mathbb{C}$ which is not on the unit circle;
	\item
There exists a prime $p\ge 2$ such that $f$ admits a root $\lambda $ in the algebraic closure $\overline{\mathbb{Q}}_p$ of the field of $p$-adic numbers $\mathbb{Q}_p$ such that $|\lambda |_p >1$.
	\end{enumerate}
	\end{lem}

	\begin{proof}
 
Suppose $r=1$, i.e., $f(z)=f_0+f_1z$. If the only root of $f$ lies on the unit circle, then necessarily $|f_0|=|f_1|$. The irreducibility of $f$ implies that $|f_0|=|f_1|=1$, and thus $\mathsf{m}(f)=0$ by \eqref{eq:mahler1}. Therefore, assuming $\mathsf{m}(f)>0$ for irreducible $f$ with $\deg(f)=1$, we conclude that $|f_0|\ne |f_1|$, and hence condition (1) must hold.

\medskip Suppose $r\ge 2$ and the roots of $f$ are all on the unit circle, i.e., suppose that condition (1) does not hold. If $|f_r|=1$ then Kronecker's theorem \cite{Kronecker} implies that all roots of $f$ are roots of unity, and \eqref{eq:mahler1} shows that $\mathsf{m}(f)=0$, in contradiction to our hypothesis. On the other hand, if $|f_r|>1$, then
	\begin{displaymath}
f(z) =f_r \bigl( z^r +\tfrac{f_{r-1}}{f_r} z^{r-1}+\cdots +\tfrac{f_0}{f_r}\bigr),
	\end{displaymath}
and Vieta's formula implies that $\bigl|\tfrac{f_0}{f_r}\bigr| = \bigl|\prod_{\zeta \in \mathbb{C}: f(\zeta )=0}\zeta \bigr| =1$, i.e., that $|f_0|=|f_r|$. Since $f$ is irreducible, $f_j/f_r$ is not integer for some $1\le j\le r-1$. Then there exists a rational prime $p$ such that $|f_j/f_r|_p>1$. Let $\lambda _i,\; 1\le i \le r$, be the roots of $f$ in the algebraic closure $\overline{\mathbb{Q}}_p$ of $\mathbb{Q}_p$, the $p$-adic rationals. By considering the $j$-th elementary symmetric polynomials of the roots and once again applying Vieta's formulae, we get that
	\begin{displaymath}
1<\biggl|\tfrac{f_j}{f_r}\biggr|_p = \biggl|\sum_{1\le i_1<i_2<\cdots <i_j \le r}\lambda _{i_1}\cdots \lambda _{i_j}\biggr|_p\le \bigl(\max_{1\le i\le r} |\lambda _i|_p\bigr)^j.
	\end{displaymath}
Thus one has $|\lambda _i|_p>1$ for some $i\in\{1,\ldots,r\}$.
	\end{proof}

The following key lemma will be used to show that for sufficiently large
$m$, the sequence of polynomials $\{z^{nm}\}_{n\ge 0}$ gives a dissociate sequence of characters of $X_f$.

	\begin{lem}[Condition (C1)]
	\label{lem:dissociate}
Suppose that $f=f_0+f_1z+\ldots +f_rz^r\in\Z[z]$ has a root in $\mathbb{C}$ or in $\overline{\mathbb{Q}}_p$ {\textup(}for some $p${\textup)} of absolute value larger than $1$. Then for any sufficiently large $m$ and any $D\ge 0$, the polynomials
	\begin{equation}
	\label{pol}
P(z) =\sum_{j=0}^D \eps_j z^{mj}, \quad \text{with }\eps_0,\eps_1,\ldots,\eps_D\in \{-2,-1,0,1,2\},
	\end{equation}
are not divisible by $f$ unless $\epsilon_0=\epsilon_1=\cdots = \epsilon_D=0$.
	\end{lem}
	\begin{proof}
First we consider the case that $f$ has a root in $\mathbb{C}$ of modulus larger than $1$. For any polynomial $g$ we introduce the notation
	\begin{displaymath}
\rho_g=\max\{|z|\mid g(z)=0\}.
	\end{displaymath}
Without loss of generality we may assume that $\eps_D\ne 0$ and consider the reduced polynomial
	\begin{displaymath}
\widetilde{P}(z) =\sum_{j=0}^D \eps_j z^{j},
	\end{displaymath}
such that $ P(z)= \widetilde{P}(z^m)$. Clearly, $\rho_{\widetilde{P}} = \rho_P^m$.

On the other hand, using the Cauchy bound on the roots of polynomials, one gets that
	\begin{displaymath}
\rho_{\widetilde{P}}\le 1+\max_{j=0,\ldots,D-1} \biggl|\frac {\eps_j}{\eps_D}\biggr|\le 3,
	\end{displaymath}
and hence $ \rho_P \le 3^\frac {1}{m}. $ Choose an integer $M\ge 1$ large enough such that $3^\frac {1}{M}<\rho_f$ (this is possible because $\rho_f >1$). Thus for all $m\ge M$, we have $\rho_P <\rho_f$. However, if $P(z)$
was divisible by $f$, we would have $\rho_f \le \rho_P$, thus arriving to a contradiction.

If $f$ has a root in $\overline{\mathbb{Q}}_p$ (for some prime $p$) of absolute value larger than $1$, the same argument works with $|\cdot|$ replaced by $|\cdot|_p$. Indeed, suppose $\zeta$ is a root of $f$ with $|\zeta|_p>1$. If $f|P$, then $\zeta$ is also a root of $P$,
so one has
	\begin{align*}
|\zeta|_{p}^{mD}&=|\zeta^{mD}|_p=\biggl| \sum_{j=0}^{D-1} \frac{\eps_j}{\eps_D} \zeta^{mj}\biggr|_p \le \max_{j=0,\ldots, D-1} \biggl| \frac{\eps_j}{\eps_D} \zeta^{mj}\biggr|_p
	\le \biggl(\max_{j=0,\ldots, D-1} \biggl| \frac{\eps_j}{\eps_D}\biggr|_p\biggr) |\zeta|_{p}^{m(D-1)}. 	\end{align*}
Thus arriving at a contradiction.
	\end{proof}

	\begin{lem}[Condition (C2)]
	\label{lem:Riesz2}
Suppose that $f(z)= f_rz^r+\ldots + f_1z+f_0 \in \mathbb{Z}[z]$ has a root in $\mathbb{C}$ or in $\overline{\mathbb{Q}}_p$ {\textup(}for some prime $p${\textup)} of absolute value larger than $1$. 		Then for all sufficiently large integers $m$, any integer $k$ with $1\le k <m$, every $D\ge0$, and all arbitrary $(D+1)$-tuples $\eps=(\eps_0,\ldots,\eps_D)$ and $\delta=(\delta_0,\ldots,\delta_D)$ in $\{-1,0,1\}^{D+1}$, the polynomial
	\begin{equation}
	\label{polQ}
Q(z) = \sum_{j=0}^D \eps_j z^{mj}-\sum_{j=0}^D \delta_j z^{mj+k}
	\end{equation}
is not divisible by $f$ unless $Q(z)\equiv 0$, i.e., unless all $\epsilon_j$'s and $\delta_j$'s are equal to zero.
	\end{lem}

	\begin{proof}
Assume $Q$ is divisible by $f$. We treat the complex case first. Namely, assume $\zeta \in \mathbb{C}$ is such that $f(\zeta)=0$ and $R\coloneqq |\zeta|>1$. Without loss of generality we may assume that $|\eps_D|+|\delta_D|>0$. We distinguish two cases.

\smallskip \emph{Case I. $\delta_D\ne 0$}. If the polynomial $Q(z)$, defined by \eqref{polQ}, is divisible by $f$, then $Q(\zeta)=0$, in other words,
	\begin{equation}
	\label{eq:lacune1}
\smash[b]{\delta_D \zeta^{mD+k} =\sum_{j=0}^D \eps_j \zeta^{mj}-\sum_{j=0}^{D-1} \delta_j \zeta^{mj+k}.}
	\end{equation}
It follows that
	\begin{displaymath}
\smash{R^{mD+k} \le \sum_{j=0}^D R^{mj}+\sum_{j=0}^{D-1} R^{mj+k}= \frac{R^{(D+1)m}-1}{R^m-1}+ \frac {R^{Dm}-1}{R^m-1}\cdot R^k,}
	\end{displaymath}
and hence
	\begin{displaymath}
\smash[t]{R^{k} < \frac{R^{m} }{R^m-1}+ \frac {R^k}{R^m-1}}.
	\end{displaymath}
As $m\to\infty $, the right hand side of this inequality converges to $1$, but the left hand side remains equal to $R^k>1$. If $m$ is large enough our assumption that $Q$ is divisible by $f$ leads to a contradiction.

\smallskip \emph{Case II. $\delta_D=0$ but $\eps_D\ne 0$.} In this case we have
	\begin{equation}
	\label{eq:lacune2}
\eps_D \zeta^{mD} =-\sum_{j=0}^{D-1}\eps_j \zeta^{jm}+\sum_{j=0}^{D-1}\delta_j \zeta^{jm+k}.
	\end{equation}
It follows that
	\begin{displaymath}
1 < \frac{1}{R^m -1} +\frac{R^k}{R^m -1}=\frac {R^k+1}{R^m-1}\le \frac{R^{m-1}+1}{R^m-1}.
	\end{displaymath}
Since $R>1$, the last inequality is violated for all sufficiently large $m$, and we again arrive at a contradiction with our assumption that $Q$ is divisible by $f$.

\smallskip In the $p$-adic case the argument is simpler because of the non-archimedean triangle inequality $|\zeta+\xi|_p\le \max(|\zeta|_p,|\xi|_p)$. Indeed, from \eqref{eq:lacune1} we get that $|\zeta|_p^{mD+k}\le |\zeta|_p^{mD}$ (impossible), and from \eqref{eq:lacune2} we get that $|\zeta|_p^{mD}\le |\zeta|_p^{m(D-1)+k}$ (equally impossible).
	\end{proof}

	\begin{proof}[Proof of Theorem \ref{thm:mainZ1}] By Remark \ref{r:irreducible} we may assume that the polynomial $f(z) = f_0+\dots +f_rz^r\in \mathbb{Z}[z]$ is primitive, irreducible, and has positive Mahler measure (cf. \eqref{eq:mahler1}). If $f$ has a root in $\mathbb{C}$ of absolute value larger than 1, the Lemmas \ref{lem:dissociate} and \ref{lem:Riesz2} show that $f$ is $m$-good for some $m\in\N$ (in fact, for all sufficiently large $m$). An application of Theorem \ref{thm:mgood} completes the proof of Theorem \ref{thm:mainZ1} in this case.

If all complex roots of $f$ have absolute value $\le1$, and if $\zeta $ is a root of $f$, we denote by $\mathbf{k}=\mathbb{Q}(\zeta )$ the algebraic number field generated by $\zeta $ and write $P_\mathbf{k}$ for the set of \emph{places} (or \emph{valuations}) of $\mathbf{k}$. If $|\zeta |_v$ is the absolute value of $\zeta $ at a place $v\in P_\mathbf{k}$, then our last assumption implies that $|\zeta |_v\le 1$ for every infinite place $v\in P_\mathbf{k}$. Since the product formula for algebraic number fields shows that $\prod_{v\in P_\mathbf{k}}|\zeta |_v=1$, we either have that $|\zeta |_v=1$ for every $v\in P_\mathbf{k}$, or that there exists a \textit{finite} place $v\in P_\mathbf{k}$ with $|\zeta |_v>1$. In the first case  $\zeta $ is a root of unity, which makes $f$ cyclotomic and $\mathsf{m}(f)=0$, contrary to our hypothesis. In the second case we can once again use the Lemmas \ref{lem:dissociate} and \ref{lem:Riesz2} to show that $f$ is $m$-good for all sufficiently large $m$, and Theorem \ref{thm:mgood} completes the proof of Theorem \ref{thm:mainZ1} as above.
	\end{proof}

\section{Bohr chaoticity of $(X_f, \alpha _f)$: the case of $d\ge 2$}
	\label{s:d>1}

This section is devoted to the proof of Theorem \ref{thm:mainZd} for $d\ge2$, which will again be based Theorem \ref{thm:mgood}. As before we assume that the polynomial $f\in R_d$ is primitive and irreducible.

\subsection{Homoclinic points of atoral polynomials in $R_d$ and the gap theorem} For every $t\in \mathbb{T}$ we set
	\begin{displaymath}
|\negthinspace|t|\negthinspace|= \min_{q\in \mathbb{Z}}|t-q|.
	\end{displaymath}

	\begin{defn}
	\label{d:homoclinic}
A point $x\in X_f$ is \textit{homoclinic} (or, more precisely, \emph{homoclinic to $0$}) if $\lim_{\bn\to\infty }|\negthinspace|x_\bn|\negthinspace|=0$. A homoclinic point $x\in X_f$ is \textit{summable} if $\sum_{\bn\in \mathbb{Z}^d}|\negthinspace|x_\bn|\negthinspace|<\infty $.
	\end{defn}

The existence of nonzero summable homoclinic points of $(X_f,\alpha _f)$ is equivalent to \emph{atorality} of the polynomial $f$:

	\begin{thm}[\cite{LSV}]
	\label{thm:atoral}
If $0\ne f\in R_d$, the following conditions are equivalent:
	\begin{enumerate}
	\item
The principal algebraic action $(X_f,\alpha _f)$ has a nonzero summable homoclinic point;
	\item
The Laurent polynomial $f$ is atoral in the sense of Definition \ref{d:atoral}.
	\end{enumerate}
	\end{thm}

For a principal algebraic $\mathbb{Z}^d$-action $(X_f,\alpha _f)$, the existence of summable homoclinic points has a number of important consequences (cf. \cite{LS}): it implies positivity of entropy and very strong specification properties of the action, and it guarantees the coincidence of entropy with the logarithmic growth rate of the number of periodic points of $\alpha _f$ (i.e., of points in $X_f$ with finite orbits under $\alpha _f$ -- cf. \cites{LS, LSV}). Somewhat surprisingly, it also plays a role in the \emph{gap theorem} stated below, which will imply the conditions (C1) and (C2) in Definition \ref{d:m-good}.

We remark in passing that some of these consequences of atorality also hold for toral polynomials, but with considerably harder proofs and/or weaker conclusions --- cf. e.g., \cite{DIM} or \cite{LSV}). However, it is not known if specification or gap properties hold in the toral case.

In order to state the gap theorem referred to above we consider, for any nonempty subset $\S\subset\Z^d$ and any integer $H\ge 1$, the set $\mathcal P(\S,H)\subset R_d$ of all Laurent polynomials with support in $\S$ and coefficients bounded in absolute value by $H$:
	\begin{displaymath}
\mathcal P(\S,H)=\bigl\{ v\in R_d\mid \textsf{supp}(v)\subseteq \S\enspace \textup{and}\enspace \|v\|_\infty \le H\}.
	\end{displaymath}
For every $\bn=(n_1,\dots ,n_d)\in \mathbb{Z}^d$ we set $\|\bn\|=\max \{|n_1|,\dots ,|n_d|\}$. Then the following is true.

	\begin{thm}[Gap Theorem]
	\label{thm:separation}
Suppose that $g\in R_d$ is primitive, irreducible and atoral. For every $H\ge1$ there exists an integer $m\ge1$ with the following property: for every pair of sets $\S,\S'\subset \mathbb{Z}^d$ with distance
	\begin{displaymath}
d(\S,\S') \coloneqq \min_{\bn\in S,\bn'\in S'}\|\bn-\bn'\|\ge m,
	\end{displaymath}
and for every $v=\sum_{\bn\in \S\cup\S'} v_{\bn}\bz^{\bn}
\in \mathcal{P}(\S\cup \S',H)$ which is divisible by $g$, the restriction of $v$ to $\S$
	\begin{equation}
	\label{eq:v_S}
v_\S=\sum_{\bn\in \S} v_{\bn}\bz^{\bn}
	\end{equation}
is also divisible by $g$.
	\end{thm}

For the proof of Theorem \ref{thm:separation} we consider the algebra $\ell ^1(\mathbb{Z}^d,\mathbb{R})$ of all {{functions}} $v\colon \bn \mapsto v_\bn$ from $\mathbb{Z}^d$ to $\mathbb{R}$ with $\|v\|_1=\sum_{\bn\in \mathbb{Z}^d}|v_\bn|<\infty $, furnished with its usual multiplication (or convolution) $(v,w)\mapsto v\cdot w$ and involution $w\mapsto w^*$, given by
	\begin{equation}
	\label{eq:convolution}
\smash[b]{(v\cdot w)_\bn = \sum_{\bm\in \mathbb{Z}^d}v_{\bm}w_{\bn - \bm} = \sum_{\bm\in \mathbb{Z}^d}v_{\bn-\bm}w_\bm,}
	\end{equation}
and
	\begin{equation}
	\label{eq:involution}
w^*_\bm = w_{-\bm}
	\end{equation}
for every $v,w\in \ell ^1(\mathbb{Z}^d,\mathbb{R})$ and $\bm,\bn\in \mathbb{Z}^d$. If we denote by $\ell ^1(\mathbb{Z}^d,\mathbb{Z})\subset \ell ^1(\mathbb{Z}^d,\mathbb{R})$ the set of all integer-valued elements of $\ell ^1(\mathbb{Z}^d,\mathbb{R})$ and identify every $h = \sum_{\bn\in \mathbb{Z}^d}h_\bn\boldsymbol{z}^\bn\in R_d$ with the element $(h_\bn)_{\bn\in \mathbb{Z}^d}\in \ell ^1(\mathbb{Z}^d,\mathbb{Z})$, we obtain an embedding
	\begin{displaymath}
R_d = \ell ^1(\mathbb{Z}^d,\mathbb{Z}) \subset \ell ^1(\mathbb{Z}^d,\mathbb{R})
	\end{displaymath}
in which the multiplication $(h,h')\mapsto h\cdot h'$ of Laurent polynomials extends to the composition \eqref{eq:convolution} in $\ell ^1(\mathbb{Z}^d,\mathbb{R})$. In fact, the multiplication $(v,w)\mapsto v\cdot w$ in \eqref{eq:convolution} is also well-defined for $w\in \ell ^1 (\mathbb{Z}^d,\mathbb{R})$ and $v\in \ell ^\infty (\mathbb{Z}^d,\mathbb{R})$, the space of all bounded sequences $(v_\bn)_{\bn\in \mathbb{Z}^d}$ in the supremum norm $\|v\|_\infty =\sup_{\bn\in \mathbb{Z}^d}|v_\bn|$, and
	\begin{displaymath}
\|v\cdot w\|_\infty \le\|v\|_\infty \|w\|_1
	\end{displaymath}
for all $w\in \ell ^1 (\mathbb{Z}^d,\mathbb{R})$ and $v\in \ell ^\infty (\mathbb{Z}^d,\mathbb{R})$.

\smallskip The shift action $\bar\sigma $ of $\mathbb{Z}^d$ on $\ell ^\infty (\mathbb{Z}^d,\mathbb{R})$, defined exactly as in \eqref{eq:shift} by
	\begin{equation}
	\label{eq:sigma2}
(\bar\sigma ^\bm v)_\bn = v_{\bm+\bn}
	\end{equation}
for every $\bm\in \mathbb{Z}^d$ and $v\in\ell ^\infty (\mathbb{Z}^d,\mathbb{R})$, extends to an action $w\mapsto w(\bar\sigma )$ of $\ell ^1(\mathbb{Z}^d,\mathbb{R})$ on $\ell^\infty (\mathbb{Z}^d,\mathbb{R})$ by bounded linear operators with
	\begin{displaymath}
w(\bar\sigma )= \sum_{\bm\in \mathbb{Z}^d}w_\bm\bar\sigma ^\bm\colon \ell^\infty (\mathbb{Z}^d,\mathbb{R})\to \ell^\infty (\mathbb{Z}^d,\mathbb{R})
	\end{displaymath}
for every $w\in \ell ^1(\mathbb{Z}^d,\mathbb{R})$. Equation \eqref{eq:sigma2} implies that
	\begin{displaymath}
\bigl(w(\bar\sigma )v\bigr)_\bn=\sum_{\bm\in \mathbb{Z}^d}w_\mathbf{m}(\bar\sigma ^\bm v)_\bn = \sum_{\bm\in \mathbb{Z}^d} w_\bm v_{\bm+\bn} = (w^*\cdot v)_\bn,
	\end{displaymath}
so that
	\begin{equation}
	\label{eq:sigma3}
w(\bar\sigma )v = w^*\cdot v
	\end{equation}
for every $w\in \ell ^1(\mathbb{Z}^d,\mathbb{R})$ and $v\in \ell ^\infty (\mathbb{Z}^d,\mathbb{R})$ (cf. \eqref{eq:involution}).

We define a surjective group homomorphism $\eta \colon \ell ^\infty (\mathbb{Z}^d,\mathbb{R}) \to \mathbb{T}^{\mathbb{Z}^d}$ by setting
	\begin{equation}
	\label{eq:eta}
\eta (v)_\bn = v_\bn \;(\textup{mod}\;1)
	\end{equation}
for every $v=(v_\bn)_{\bn\in \mathbb{Z}^d}$ and $\bn\in \mathbb{Z}^d$. Note that $\eta $ is shift-equivariant in the sense that
	\begin{displaymath}
\eta \circ \bar\sigma ^\bn = \sigma ^\bn \circ \eta
	\end{displaymath}
for every $\bn\in \mathbb{Z}^d$; more generally, if $w\in R_d=\ell ^1(\mathbb{Z}^d,\mathbb{Z})$, then
	\begin{equation}
	\label{eq:eta-f}
\eta \circ w(\bar\sigma ) = w(\sigma )\eta.
	\end{equation}

\smallskip For every $x\in \mathbb{T}^{\mathbb{Z}^d}$ there exists a unique point $x^\# \in (-\frac12,\frac12]^{\mathbb{Z}^d}\subset \ell ^\infty (\mathbb{Z}^d,\mathbb{R})$, called the \emph{lift} of $x$, such that
	\begin{equation}
	\label{eq:lift}
\eta (x^\#) = x.
	\end{equation}

Let $g=\sum_{\bn\in \mathbb{Z}^d} g_\bn \boldsymbol{z}^\bn \in R_d$ be the Laurent polynomial appearing in the statement of Theorem \ref{thm:separation} and set
	\begin{displaymath}
f=g^* = \sum_{\bn\in \mathbb{Z}^d} g_\bn \boldsymbol{z}^{-\bn }.
	\end{displaymath}
Since $U(f)= U(g)$, $f$ is again atoral and has nontrivial summable homoclinic points by Theorem \ref{thm:atoral}.

	\begin{lem}
	\label{lem:lift}
For every $x\in \mathbb{T}^{\mathbb{Z}^d}$, the following is true:
	\begin{enumerate}
	\item
$x\in X_f$ if and only if $f(\bar{\sigma })x^\# \in \ell ^\infty (\mathbb{Z}^d,\mathbb{Z}),\ i.e. \ x^\#\cdot f^* \in \ell ^\infty (\mathbb{Z}^d,\mathbb{Z})
$
{\upshape(}cf. \eqref{eq:sigma3}{\upshape)};
	\item
$x$ is a nontrivial summable homoclinic point of $\alpha _f$ if and only if $x^\#\in \ell ^1(\mathbb{Z}^d,\mathbb{R})$,
$h\coloneqq x^\#\cdot f^* \in \ell ^1(\mathbb{Z}^d,\mathbb{Z})=R_d$ and $h$ is not divisible by $f^*$ in $R_d$.
	\end{enumerate}
	\end{lem}

	\begin{proof}
(1) Suppose that $x\in \T^{\Z^d}$. By \eqref{eq:eta-f}, we have
	\begin{displaymath}
\eta (f(\bar{\sigma })x^\#)=f(\sigma )\eta (x^\#)=f(\sigma )x.
	\end{displaymath}
So, $x\in X_f$, i.e. $f(\sigma )x=0$ if and only if $f(\bar{\sigma })x^\# \in \ell ^\infty (\mathbb{Z}^d,\mathbb{Z})$.

(2) If $x$ is a nontrivial summable homoclinic point of $\alpha _f$, then $x^\# \in \ell ^1(\mathbb{Z}^d,\mathbb{R})$, and part (1) of this proof implies that $h=x^\#\cdot f^*\in \ell ^1(\mathbb{Z}^d,\mathbb{Z}) = R_d$. If $h$ were divisible by $f^*$, that is, if $x^\#\cdot f^* = h\cdot f^*$ for some $h\in R_d$, then $(x^\#-h)\cdot f^*=0$. Then \cite{LP}*{Theorem 2.1} would imply that $x^\#=h$ and hence $x=\eta (x^\#)=\eta (h)=0$, contrary to our conditions on $x$. The converse is obvious.
	\end{proof}

	\begin{proof}[Proof of Theorem \ref{thm:separation}]
Since $f=g^*$ is atoral, there exists a nontrivial summable homoclinic point $x\in X_f$. Let $x^\#\in \ell ^1(\mathbb{Z}^d,\mathbb{R})$ be the lift of $x$ (cf. \eqref{eq:lift}), and let $h=f(\bar{\sigma })x^\#=x^\#\cdot f^* \in R_d$ (cf. Lemma \ref{lem:lift}). Since $x^\#\in \ell ^1(\mathbb{Z}^d,\mathbb{R})$, there exists an integer $R=R(x,f,H)$ such that
	\begin{displaymath}
\smash[t]{\sum_{\|\bn\|\ge R} |x^\#_{\bn}| <\frac 1{2H\|f\|_1}.}
	\end{displaymath}
For every nonempty subset $\S\subset \mathbb{Z}^d$ we set
	\begin{displaymath}
B_R(\S)=\{\bn\in \mathbb{Z}^d\mid d(\bn,\S)=\min_{\bn'\in \S}\|\bn-\bn'\|\le R\}.
	\end{displaymath}

Let $\S,\S'\subset \mathbb{Z}^d$ be two subsets of $\mathbb{Z}^d$ with distance $d(\S,\S')\ge 3R$. Suppose that a Laurent polynomial $v\in\mathcal P(\S\cup\S', H)$ is divisible by $f^*$, i.e., that $v=\phi \cdot f^*$ for some $\phi \in R_d$. Then
	\begin{enumerate}
	\item[(i)]
$v\cdot x^\#\in R_d$;
	\item[(ii)]
$\textsf{supp}(v\cdot x^\#)\subset B_R(\S)\cup B_R(\S')$.
	\end{enumerate}
Indeed, (i) follows from Lemma \ref{lem:lift} (2):
	\begin{equation*}
	\label{eq:vw_in_Rd}
v\cdot x^\# = (\phi \cdot f^*) \cdot x^\#= \phi\cdot (f^*\cdot x^\#)=\phi\cdot h\in R_d,
	\end{equation*}
because both $\phi $ and $h$ belong to $R_d$. For (ii) we note that every $\bn \notin B_R(\S)\cup B_R(\S')$ satisfies that $d(\bn, \S\cup\S')> R$. Then $v_{\bn-\bm}=0$ for all $\bm$ with $\|\bm\|\le R$, and hence
	\begin{displaymath}
|(v\cdot x^\#)_{\bn}|=\biggl|\sum_{\bm\in\Z^d} x^\#_{\bm}v_{\bn-\bm}\biggr|\le \|v\|_\infty \sum_{\|\bm\|>R}| x^\#_{\bm}| < H\cdot \frac {1}{2H\|f\|_1}\le \frac 12.
	\end{displaymath}
Since $(v\cdot x^\#)_{\bn} \in \mathbb{Z}$ by (i), it follows that $(v\cdot x^\#)_{\bn}=0$.

Let $\psi $ be the restriction of $v \cdot x^\#$ to $B_R(\S)$, and let $v_\S$ and $v_{\S'}$ be the restrictions of $v$ to $\S$ and $\S'$, respectively. Then $\psi \in R_d$ by (i), and we claim that
	\begin{equation}
\|\psi - v_\S\cdot x^\#\|_{\infty} < \frac{1}{2\|f\|_1},
	\end{equation}
i.e., that
	\begin{equation}
	\label{eq:vs2}
|\psi_\bn- (v_\S\cdot x^\#)_\bn| < \frac{1}{2\|f\|_1} \quad \textup{for every}\enspace \bn \in \Z^d.
	\end{equation}

\smallskip Indeed, if $\bn \in B_R(\S)$, then $d(\bn,B_R(\S')) \ge R$, and hence
	\begin{equation}
	\label{eq:vs3}
|(v_{\S'}\cdot x^\#)_\bn| = \biggl|\sum_{\bm\in \S'} v_\bm x^\#_{\bn -\bm}\biggr|\le \|v\|_\infty \sum_{\|\bl\|\ge R} |x^\#_{\bl}|\le H\cdot \frac{1}{2 H \|f\|_1} = \frac{1}{2\|f\|_1}.
	\end{equation}
Since $\psi _\bn = (v_{\S}\cdot x^\#)_\bn + (v_{\S'}\cdot x^\#)_\bn$, it follows that
	\begin{displaymath}
|\psi _\bn - (v_{\S}\cdot x^\#)_\bn| = |(v_{\S'}\cdot x^\#)_\bn| < \frac{1}{2\|f\|_1}
	\end{displaymath}
by \eqref{eq:vs3}. On the other hand, if $\bn\notin B_R(\S)$, then
	\begin{displaymath}
|\psi _\bn - (v_\S \cdot x^\#)_\bn |=|(v_\S \cdot x^\#)_\bn | = \biggl|\sum_{\bm\in \S}(v_\bm x^\#_{\bn-\bm})\biggr| \le H\cdot \frac{1}{2 H \|f\|_1} = \frac{1}{2\|f\|_1}.
	\end{displaymath}
This proves \eqref{eq:vs2} for every $\bn\in \mathbb{Z}^d$.

Since both $v_\S\cdot x^\#\cdot f^* = v_\S\cdot h$ and $\psi $ lie in $R_d$ we have that $(\psi -v_\S\cdot x^\#)\cdot f^* \in R_d$, but the smallness of the coordinates of $\psi -v_\S\cdot x^\#$ in \eqref{eq:vs2} implies that $(\psi -v_\S\cdot x^\#)\cdot f^*=0$. Thus we have proved that $\psi \cdot f^*=v_\S\cdot x^\#\cdot f^*=v_\S\cdot h$, where $h$ is not divisible by $f^*$ (cf. Lemma \ref{lem:lift} (2)). As $g=f^*$ is irreducible, we have proved that $v_\S$ is divisible by $g$, as claimed in the statement of this theorem.

This completes the proof of Theorem \ref{thm:separation} with $m\ge 3R$.
	\end{proof}

\subsection{The conditions (C1) and (C2): divisibility by $f$ of lacunary polynomials}

	According to Theorem \ref{thm:mgood}, in order to prove Theorem \ref{thm:mainZd}, it suffices to prove that any irreducible and atoral polynomial $f\in R_d$ is $m$-good for a sufficiently large $m\in\N$. Now we are going to prove this and finish the proof of Theorem \ref{thm:mainZd}.

Theorem \ref{thm:separation} has an immediate corollary which implies that any atoral polynomial is $m$-good for sufficiently large $m$.

	\begin{cor}
	\label{cor:separation1}
Suppose that $f\in R_d$ is irreducible and atoral, and that $|\textsf{supp}(f)| > 1$. Then there exists, for every $H\ge 1$, an integer $m\ge1$ with the following property: for any set $\S\in\Z^d$ which is $m$-separated in the sense that
	\begin{displaymath}
\|\bk-\bn\|\ge m
\text{ for any pair }\bk, \bn\in\S,\ \bk\ne\bn,
	\end{displaymath}
no non-zero polynomial $g\in\mathcal P(\S,H)$ is divisible by $f$.
	\end{cor}

	\begin{proof}
For $H\ge 1$ and $f$ fixed, choose $m$ as in the statement of Theorem \ref{thm:separation} (i.e. $m\ge 3R$ in the proof of that theorem). Consider an arbitrary $m$-separated set $\S$ and any
non-trivial polynomial $v = \sum_{\bn\in \S}v_\bn\boldsymbol{z}^\bn\in\mathcal P(\S,H)$.

If $|\textsf{supp}(v)|=1$, then $v$ cannot be divisible by $f$, since $|\textsf{supp}(f)|>1$ by assumption. Assume therefore that $|\textsf{supp}(v)|\ge 2$, and that $v$ is divisible by $f$. Since for any $\bn\in \textsf{supp}(v)$, the sets
	\begin{displaymath}
\mathcal{T}=\{\bn\}, \qquad \mathcal{T}'=\textsf{supp}(v)\smallsetminus \{\bn\}
	\end{displaymath}
have distance at least $m$ and hence, by Theorem \ref{thm:separation}, the restriction of $v$ to $\mathcal{T}$, i.e. $v_{\mathcal{T}}=v_{\bn}\bz^{\bn}$
must be divisible by $f$, which is impossible. So, $v$ is not divisible by $f$.
	\end{proof}

The condition that $|\textsf{supp}(f)|>1$ in Corollary \ref{cor:separation1} is obviously necessary: the polynomial $f=2$ is obviously irreducible and atoral, and divides $2g$ for every $g\in R_d$ (irrespective of whether $g$ is $m$-separated or not).

	\begin{cor}
	\label{cor:separation2}
Suppose that $f\in R_d$ is irreducible and atoral, and that $|\textsf{supp}(f)| > 1$. For all sufficiently large $m\ge 1$ and every $\boldsymbol{k}\in [0,m-1]^d\smallsetminus\{\bo\}$, no $v\in \mathcal{P}\bigl(m\mathbb{Z}^d \cup (m\mathbb{Z}^d+{\boldsymbol k}),1\bigr)$ with $v \ne 0$ is divisible by $f$.
	\end{cor}

	\begin{proof}
Put $H=1$ and let $m\ge 6R$, where $R$ is the number appearing in the proof of Theorem \ref{thm:separation}.
Suppose $v\in \mathcal{P}\bigl(m\mathbb{Z}^d \cup (m\mathbb{Z}^d+\boldsymbol{k}),1\bigr)$ is a non-trivial polynomial divisible by $f$. Consider the decomposition
$\textsf{supp}(v)=\S_0\sqcup\S_1$ where
	\begin{displaymath}
\S_0=\textsf{supp}(v)\cap m\Z^d, \quad \S_1=\textsf{supp}(v)\cap
(m\Z^d+\bk).
	\end{displaymath}
Both sets $\S_0,\S_1$ are $m$-separated, as subsets of $m\Z^d$ and $m\Z^d+\bk$ respectively.

We claim that for any $\bn\in \S_0$ there exists $\bn'=\bn'(\bn)\in \S_1$ such that $d(\bn,\bn')<3R$. Otherwise, there exists $\bn\in \S_0$ such that $d(\bn, \S_1)\ge 3R$ so that $d(\bn, \textsf{supp}(v)\smallsetminus \{\bn\})\ge 3R$. Then, by Theorem \ref{thm:separation}, the restriction of $v$ to $\{\bn\}$, i.e., $\pm \bz^{\bn}$, is divisible by $f$, which is impossible. Similarly, for any $\bn'\in\S_1$ there exists $\bn\in S_0$ such that $d(\bn,\bn')<3R$. Thus the support of $v$ is a union of distinct pairs:
	\begin{displaymath}
\textsf{supp}(v)=\bigcup_{\bn\in\S_0}\{\bn,\bn'\},
	\end{displaymath}
where the distance within each pair is at most $3R$.

Given a pair $\{\bn,\bn'\}$, consider the decomposition of $\textsf{supp}(v)$:
	\begin{displaymath}
\S=\{\bn,\bn'\}, \quad \S'=\textsf{supp}(v)\smallsetminus \S.
	\end{displaymath}
The fact that $m\ge 6R$ implies $d(\S, \S')\ge 3R$. Indeed, $d(\bn, \S')= d(\bn, \bn^*)$ for some $\bn^*\in \S'$ and
	\begin{eqnarray*}
d(\bn, \bn^*) & \ge & m \ \ \mbox{\rm if}\ \bn^*\in \S_0;
	\\
d(\bn, \bn^*)& \ge & d(\bn', \bn^*) - d(\bn', \bn) \ge m -3R \ \ \mbox{\rm if}\ \bn^*\in \S_1.
	\end{eqnarray*}
It follows that $d(\bn, \S')>3R$. Similarly, $d(\bn', \S')>3R$.

Applying Theorem \ref{thm:separation} to $\S$ and $\S'$, we conclude that the restriction of $v$ to $\S=\{\bn,\bn'\}$, i.e.
	\begin{displaymath}
v_{\S} =v_{\bn} \bz^{\bn}+v_{\bn'} \bz^{\bn'}, \quad v_{\bn},v_{\bn'}\in\{-1,1\},
	\end{displaymath}
must be divisible by $f$, which is impossible, since $v_\S$ is of the form
	\begin{displaymath}
\pm \bz^{\bm}(1\pm \bz^{\bl}), \quad \bm\in\Z^d,\ \bl\in\mathbb{N}^d,
	\end{displaymath}
and hence is a product of a unit ($\pm \bz^{\bm}$) and a generalized cyclotomic polynomial $(1\pm \bz^{\bl})$, and thus must have zero Mahler measure $m(v_S)=0$. This implies that $\mathsf{m}(f) =0$, in violation of Theorem \ref{thm:atoral}.
	\end{proof}

	\begin{proof}[Proof of Theorem \ref{thm:mainZd}] Since Remark \ref{r:irreducible} allows us to assume without loss in generality that the polynomial $f\in R_d$ is primitive and irreducible, the proof of Bohr chaoticity under the additional assumption of atorality of $f$ is now complete.

Indeed, if $|\textsf{supp}(f)|=1$, atorality implies that we are in the situation of Example \ref{ex:constant} with $p>1$, so that $(X_f,\alpha _f)$ is Bohr chaotic. If $|\textsf{supp}(f)|\ge 2$, Corollary \ref{cor:separation1} for $H=2$ and Corollary \ref{cor:separation2} show that the conditions (C1) and (C2) are satisfied. Therefore, Bohr chaoticity of $(X_f,\alpha _f)$ for irreducible atoral polynomials $f\in R_d$ follows from Theorem \ref{thm:mgood}.
	\end{proof}

\section{Concluding remarks: toral polynomials}\label{s:atoral}

We have shown that a principal $\Z$-action is Bohr chaotic if and only if
it has positive entropy. A principal $\Z^d$-action, $d>1$, {{is}} shown to be Bohr chaotic
if it has positive entropy and is atoral.
We believe that the atorality assumption (equivalently, existence of a non-trivial summable homoclinic point)
can be removed.

Irreducible toral polynomials come in two flavours: those, for which $X_f$ has no non-trivial homoclinic points (cf. e.g., \cite{LS}*{Example 7.1}), and those, for which $X_f$ has no summable homoclinic points, but uncountably many nonzero homoclinic points $v\in X_f$ with the property that $v^\# \cdot f^* \in f^* R_d$ (cf. e.g., \cite{LS}*{Example 7.3}; for notation we refer to \eqref{eq:lift}). Unfortunately, none of these latter homoclinic points can be used in our proof of the Gap Theorem (Theorem \ref{thm:separation}), since the key Lemma \ref{lem:lift} is not valid in this case.

Remarkably, for toral examples of the kind illustrated in \cite{LS}*{Example 7.1} we can still prove Bohr chaoticity by using the Theorems \ref{thm:mgood} and \ref{thm:mgood_d=1}.

	\begin{thm}
	\label{thm:toral}
Let $g\in R_1$ be an irreducible polynomial with positive Mahler measure and with all roots of absolute value $1$, and define $f\in R_d$ by $f(z_1,\ldots,z_d)=g(z_1)$. Then $(X_f,\alpha _f)$ is Bohr chaotic.
	\end{thm}

	\begin{rem}
	\label{r:toral}
There exist infinitely many distinct irreducible polynomials $g\in R_1$ with the properties required in Theorem \ref{thm:toral}. In order to see this we follow a short note from \texttt{math.stackexchange.com} \cite{stack}:

Let $\vartheta $ be a totally real algebraic number all of whose conjugates $\vartheta _1=\vartheta ,\vartheta _2,\dots ,\vartheta _r$ have absolute values strictly less than 2 (to find such a $\vartheta $, take any totally real algebraic number and divide by a big integer). Assume also that $\vartheta $ is not itself an algebraic integer. Let $\beta $ be a solution to the equation
	\begin{displaymath}
\beta +1/\beta =\vartheta .
	\end{displaymath}
All the conjugates of $\vartheta $ are, by assumption, real numbers in the interval $(-2,2)$. This forces all the conjugates of $\beta $ to be complex numbers of absolute value one --- which is what we want. Moreover, $\beta $ will not be a root of unity, since otherwise $\vartheta $ would be an algebraic integer. Then $\beta $ is a root of the polynomial
	\begin{displaymath}
h(x)=\prod _{i=1}^r (x^2 - \vartheta _ix + 1),
	\end{displaymath}
Clearing denominators in $h$, one gets the desired polynomial $g$. It will be irreducible, because, by looking at infinite places, $[\mathbb{Q}(\beta ):\mathbb{Q}(\vartheta )]=2$.

For example, $\vartheta =1/2$ yields $g=2z^2-1+2$, $\vartheta = -6/5$ yields $g=5z^2-6z+5$, $\vartheta = 1/\sqrt 2$ yields $g=2z^4+3z^2+2$, etc.
	\end{rem}

	\begin{proof}[Proof of Theorem \ref{thm:toral}]
Since $g$ with $\mathsf m(g)>0$ is $m$-good for some sufficiently large $m$ by Theorem \ref{thm:mgood_d=1}, and $f(z_1,\ldots,z_d)=g(z_1)$ is also $m$-good, but now viewed as a polynomial in $d$-variables. Indeed, if $f(\bz)=g(z_1)$ divides a non-trivial polynomial $h$ of the form
	\begin{align*}
h(\bz) &=\sum_{\bn \in \mathbb{Z}^d} \varepsilon _\bn \bz^{m\bn}
	\end{align*}
with $\varepsilon _\bn=\varepsilon_{(n_1,n_2,\ldots,n_d)} \in \{-2,-1,0,1,2\}$, then
by rewriting $h$ as
	\begin{align*}
h(\bz)
&= \sum_{(n_2,\ldots,n_{d})\in \Z^{d-1}}
\left( \sum_{n_1\in\Z} \varepsilon_{(n_1,n_2,\ldots,n_d)} z_1^{m n_1}\right) z_2^{m n_2}\cdots z_{d}^{m n_d}
	\\
&=: \sum_{(n_2,\ldots,n_{d})\in \Z^{d-1}}
h_{n_2,\ldots,n_d}(z_1^m) z_2^{mn_2}\cdots z_{d}^{mn_d}
	\end{align*}
we conclude that $g(z_1)$ must
divide all polynomials $h_{n_2,\ldots,n_d}(z_1^m)$, some of which are nonzero. Since $g(z_1)$ is $m$-good, the resulting contradiction proves that condition (C1) is valid for $f$. For the proof of condition (C2) we can proceed similarly.

Theorem \ref{thm:toral} now follows from Theorems \ref{thm:mgood}.
	\end{proof}

Theorem \ref{thm:toral} obviously applies only to an extremely restricted class of toral polynomials. A more interesting example of an irreducible toral polynomial with positive entropy in \cite{LS}*{Example 7.3} is given by the polynomial
	\begin{displaymath}
f(z_1,z_2)= 3-z_1-\frac 1{z_1}-z_2-\frac 1{z_2}.
	\end{displaymath}
The unitary variety of $f$ is a smooth real-analytic curve
	\begin{displaymath}
\mathsf U(f)=\left\{(e^{2\pi i s},e^{2\pi i t}): t=\pm \frac{1}{2 \pi} \cos ^{-1}\left(\frac{3}{2}-\cos 2 \pi s\right), \quad-\frac{1}{6} \leq s \leq \frac{1}{6}
\right\}.
	\end{displaymath}
Moreover, $\mathsf U(f)$ is connected and has curvature bounded away from zero. As explained in \cite{LS}*{Example 7.3.}, one can use this fact to prove the existence of probability measures supported on $\mathsf U(f)$ whose Fourier transforms (as functions on $\widehat{\mathbb{S}^d}=\mathbb{Z}^d$) vanish at infinity. By translating this information back to the system $(X_f,\alpha _f)$ one obtains uncountably many homoclinic points $x\in X_f$ satisfying
	\begin{displaymath}
|x_{\bn}| \le \frac {C}{1+\|\bn\|^{\frac 12}}\text{ for all \/}\bn\in \Z^2.
	\end{displaymath}
Unfortunately, none of these homoclinic points can be used in our proof of the Gap Theorem (Theorem
\ref{thm:separation}), since they are not summable.

It would be interesting to see whether one can prove the Gap Theorem for $f=3-z_1-\frac 1{z_1}-z_2-\frac 1{z_2}$ directly, using some elementary methods, or establish Bohr chaoticity of $X_f$ by
some other means.

	\end{document}